\documentclass[12pt]{article}
\usepackage{latexsym}
\usepackage{amssymb}
\usepackage{euscript}
\let\script\EuScript

\newtheorem{theorem}{Theorem}[section]
\newtheorem{definition}[theorem]{Definition}
\newtheorem{lemma}[theorem]{Lemma}
\newtheorem{corollary}[theorem]{Corollary}
\newtheorem{notation}[theorem]{Notation}
\newtheorem{claim}[theorem]{Claim}
\def\restrict{\upharpoonright}
\def\some{\exists}
\def\comp{\circ}
\def\N{\mathbb N}
\def\R{\mathbb R}
\def\forces{\Vdash}
\def\nin{\notin}
\def\all{\forall}

\def\o1{\omega_1}
\def\la{\langle}
\def\ra{\rangle}

\def\romu{{\mathbb R}^{\omega}/U}

\def\Q{\mathbb{Q}}
\def\halmos{{}{\hfill $\Box$} \vskip 5pt \par}
\def\soo{S_{\omega_1}}
\def\goo{G_{\omega_1}}
\def\ea{\mathcal{E}}
\def\morass{\langle (\theta_\alpha \mid \alpha\leq
\omega_1),(\mathcal{F}_{\alpha
\beta}\mid \alpha<\beta\leq
\omega_1)\rangle}
\def\gaptwo{\langle \overrightarrow{\varphi}, \overrightarrow{\script{G}}, \overrightarrow{\theta}, \overrightarrow{\script{F}} \rangle}
\def\gapone{\langle \overrightarrow{\varphi}, \overrightarrow{\script{G}} \rangle}
\begin{document}
\setlength{\baselineskip}{21pt}
\title{Discontinuous Homomorphisms
of $C(X)$ with $2^{\aleph_0}>\aleph_2$}
\author{Bob A. Dumas\\
Department of Philosophy\\
University of Washington\\
Seattle, Washington 98195}
\date{November 11, 2016}

\bibliographystyle{plain}

\maketitle
\begin{abstract}
Assume that $M$ is a c.t.m.
of $ZFC+CH$ containing a simplified
$(\omega_1,2)$-morass, $P\in M$ is the
poset adding $\aleph_3$ generic reals
and $G$ is $P$-generic over $M$.
In $M$ we construct a function between sets
of terms in the forcing language,
that interpreted in $M[G]$ is
an $\R$-linear
order-preserving monomorphism from the
finite elements of an ultrapower of the
reals, over a non-principal ultrafilter
on $\omega$, into the Esterle algebra
of formal power series.  Therefore
it is consistent that $2^{\aleph_0}=\aleph_3$
and, for any infinite compact
Hausdorff space $X$,
there exists a discontinuous
homomorphism of $C(X)$,
the algebra of continuous
real-valued functions on $X$.
\end{abstract}

\baselineskip = 18pt
\setcounter{section}{0}

\section{Introduction}
This paper addresses a problem in the
theory of Banach Algebras concerning
the existence of discontinuous
homomorphisms of $C(X)$, the algebra
of continuous
real-valued functions with domain
$X$, where $X$ is an infinite compact
Hausdorff space. In [\ref{Johnson}], B. Johnson
proved that there is a discontinuous
homomorphism of $C(X)$ provided that
there is a nontrivial submultiplicative
norm on the finite elements of an
ultrapower of $\R$ over $\omega$. In
[\ref{Esterle1}], J. Esterle constructs
an algebra of formal power series,
$\ea$, and shows in [\ref{Esterle2}]
that the infinitesimal elements of
$\ea$ admit a nontrivial
submultiplicative norm.  By results of
Esterle in [\ref{Esterle3}], it is
known that $\ea$ is an
$\eta_1$-ordering of cardinality
$2^{\aleph_0}$. Furthermore, $\ea$ is a
totally ordered field by a result of
Hahn in 1907 [\ref{Hahn}], and is
real-closed by a result of Maclane
[\ref{Maclane}].

It is a theorem of P. Erd\"{o}s, L.
Gillman and M. Henriksen in
[\ref{Erdos}] that any pair of
$\eta_1$-ordered real-closed fields of
cardinality $\aleph_1$ are
isomorphic as ordered fields.  In fact, it is shown
using a back-and-forth argument that
any order-preserving field isomorphism between countable
subsets of $\eta_1$-ordered real-closed
fields may be extended to an
order-isomorphism.  It is a standard
result of model theory that for any
non-principal ultrafilter $U$ on $\omega$, $\romu$ is
an $\aleph_1$-saturated real-closed
field (and hence an $\eta_1$-ordering).
By a result of Johnson [\ref{Johnson}],
between any pair of $\eta_1$-ordered
real-closed fields of cardinality
$\aleph_1$ there is an $\R$-linear
order-preserving field isomorphism
(hereafter $\R$-isomorphism).
This implies, in a model of the
continuum hypothesis (CH), that there
is an $\R$-linear order-preserving monomorphism
(hereafter $\R$-monomorphism)
from the finite elements of $\romu$ to $\ea$, and hence
in models of ZFC+CH there exists a
discontinuous homomorphism of $C(X)$.
The proof that in a model of ZFC+CH
there exists a discontinuous
homomorphism of $C(X)$ is due
independently to Dales [\ref{Dales1}]
and Esterle [\ref{Esterle2}].

Shortly thereafter R. Solovay found a
model of ZFC+$\neg$CH in which all
homomorphisms of $C(X)$ are continuous.
Later, in his Ph.D. thesis, W.H. Woodin
constructed a model of ZFC+Martin's Axiom in
which all homomorphisms of $C(X)$ are
continuous. This
naturally gave rise to the question of
whether there is a model of set theory
in which CH fails and there is a
discontinuous homomorphism of $C(X)$.
Woodin subsequently showed that in the
Cohen extension of a model of
ZFC+CH by generic reals
indexed by $\omega_2$, there is
a discontinuous homomorphism of
$C(X)$ [\ref{Woodin}]. Woodin shows that
in this model the gaps in $\ea$ that must be witnessed in
a classical back-and-forth construction
are always countable.  He observes
that this construction may not be extended to a
Cohen extension by more than $\aleph_2$ generic
reals.  He suggests the plausibility
of using morasses to construct
an $\R$-monomorphism
from the finite elements of an
ultrapower of the reals to the
Esterle algebra in
generic extensions with more than $\aleph_2$ generic
reals.  Woodin's
argument does not extend to higher
powers of the continuum and leaves open
the question of whether there exists a
discontinuous homomorphism of $C(X)$ in
models of set theory in which
$2^{\aleph_0}>\aleph_2$.
In this paper we show that the existence of
a simplified $(\omega_1,2)$-morass in a
model of $ZFC + CH$ is sufficient for the existence
of a discontinuous homomorphism of $C(X)$
in a model in which $2^{\aleph_0}=\aleph_3$.

We show that in the Cohen extension adding
$\aleph_2$ generic reals to
a model of $ZFC + CH$ containing a simplified
$(\omega_1,1)$-morass, there is a level,
morass-commutative term function that, interpreted
in the Cohen extension, is an $\R$-monomorphism
of the finite elements of an ultrapower of $\R$
over $\omega$ into the Esterle Algebra.
This is achieved with a
transfinite construction of length $\omega_1$, utilizing the
morass functions from the gap-one morass to complete the
construction of size $\aleph_2$ by commutativity with morass maps.
Using the techniques of this argument, we construct a term function
with a transfinite construction of length $\omega_1$ and utilize morass-commutativity
with the embeddings of a gap-2 morass to complete the construction
of an $\R$-monomorphism from the finite elements
of a standard ultrapower of $\R$ over
$\omega$ to the Esterle Algebra in the Cohen extension adding
$\aleph_3$ generic reals.

The technical obstacles
to such a construction may be reduced to conditions
we call morass-extendability.  This paper is dependent on the
results of [\ref{Dumas}] and [\ref{Dumas2}], in
which term functions are constructed that are forced to be
order-preserving functions.  The arguments here are very similar,
with the additional requirement that the functions
are also ring homomorphisms.

\section{Preliminaries}
In our initial construction we use a simplified
$(\omega_1,1)$-morass.  We
construct a function on terms in the forcing
language for adding $\aleph_2$ generic
reals that is forced in all generic
extensions to be an $\R$-monomorphism from the finite
elements of $\romu$, where $U$ is a
standard non-principal ultrafilter
(see Definition 6.14 [\ref{Dumas}]) in
the generic extension, into
the Esterle Algebra, $\ea$.  In
some sense we follow the classical
route to such constructions - extension by
transcendental elements in an inductive
construction of length
$\omega_1$. We will
require commutativity with morass maps
to construct a function on a domain of
cardinality $\aleph_2$ making only
$\aleph_1$ many explicit commitments.
However with each
commitment of the construction, there are
uncountably many
future commitments implied by commutativity
with morass maps.

In [\ref{Velleman}] D. Velleman defines
a simplified $(\omega_1,1)$-morass.
\begin{definition}$($Simplified $(\omega_1,1)$-morass$)$\label{thm2.1}
A simplified $(\omega,1)$-morass is a
structure
\[ \mathcal{M}=\langle (\theta_\alpha \mid \alpha\leq
\omega_1),(\mathcal{F}_{\alpha
\beta}\mid \alpha<\beta\leq
\omega_1)\rangle \] that satisfies the following conditions:\\
(P0) (a) $\theta_0=1$,
$\theta_{\omega_1}=\omega_2$, $(\forall
\alpha<\omega_1)\  0<\theta_\alpha<\omega_1$.\\
\  (b) $\mathcal{F}_{\alpha \beta}$ is
a set of order-preserving
functions $f:\theta_\alpha \to \theta_\beta$.\\
(P1) $\mid \mathcal{F_{\alpha \beta}}
\mid \leq \omega$ for all
$\alpha<\beta<\omega_1$.\\
(P2) If $\alpha<\beta<\gamma$, then
$\mathcal{F}_{\alpha \gamma}=\{ f\circ
g\mid f\in \mathcal{F}_{\beta \gamma}$,
$g\in
\mathcal{F}_{\alpha \beta} \}$.\\
(P3) If $\alpha <\omega_1$, then
$\mathcal{F}_{\alpha (\alpha+1)}=\{
\textnormal{id}\restrict
\theta_{\alpha} , f_{\alpha} \}$ where
$f_{\alpha}$ satisfies: \[ (\some
\delta_{\alpha}<\theta_\alpha)\
f_\alpha \restrict \delta_{\alpha} =
\textnormal{id}\restrict
\delta_{\alpha} \  \textnormal{ and } \
f_\alpha(\delta_{\alpha})\geq
\theta_\alpha. \] (P4) If $\alpha \leq
\omega_1$ is a limit ordinal, $\beta_1,
\beta_2<\alpha$, $f_1\in
\mathcal{F}_{\beta_1 \alpha}$ and
$f_2\in \mathcal{F}_{\beta_2 \alpha}$,
then there is $\gamma<\alpha$,
$\gamma>\beta_1, \beta_2$, and there is
$f_1'\in \mathcal{F}_{\beta_1 \gamma}$,
$f_2'\in \mathcal{F}_{\beta_2 \gamma}$,
$g\in \mathcal{F}_{\gamma \alpha}$
such that $f_1=g\comp f_1'$ and $f_2=g\comp f_2'$.\\
(P5) For all $\alpha>0$,
$\theta_{\alpha}=\bigcup \{
f[\theta_{\beta}] \mid \beta<\alpha$,
$f\in \mathcal{F}_{\beta \alpha} \}$.
\end{definition}
Simplified gap-1 morasses, as well as higher gap
simplified morasses, are known to exist in $L$.

We will construct, by an inductive
argument of length $\omega_1$, a function
between sets of
terms in the forcing language adding
$\aleph_2$ generic reals.  We interpret
the morass functions on ordinals as
functions between terms in the
forcing language and require that the
set of terms under construction satisfy
certain commutativity constraints with
the morass functions. It is implicit
that any commitment to an ordered pair
of terms in the construction is
\emph{de facto} a commitment to
uncountably many ordered
pairs in mutually generic extensions.
In [\ref{Dumas}] we worked explicitly
with terms in the forcing language.  We
wish to simplify the details of the
construction by working with objects in
a forcing extension. We use the notions
of term complexity and strict level,
from [\ref{Dumas}], and apply it to
objects in a forcing extension.
\begin{notation}$($P(A)$)$\label{thm4.7}
If $A$ is a set of ordinals, we
let $P(A)$ be the poset adding generic
reals indexed by the ordinals of $A$.
That is,
\[ P(A):= Fn(A\times \omega,2), \]
the finite partial functions from
$A\times \omega$ to $2$.
\end{notation}
Let $M$ be a c.t.m of
ZFC, $\beta$ be an ordinal
and $P$ be the poset adding
generic reals indexed by $\beta$,
then $P(\beta)=Fn(\beta \times \omega,2)$.
Suppose $S\subseteq \beta$.
Let $\tau \in M^{P(S)}$ be a discerning term
in the forcing language adding
generic reals indexed by $S$.
Then $\tau$ has strict level $S$
provided that, for any proper
subset of $S$, $A$, and term,
$\sigma$, in the forcing language
adding generic reals indexed by
$A$, $\forces \tau \neq \sigma$.
Not all terms have strict levels.

Alternatively, we consider $P(\beta)$
as the product forcing $P(A)\times
P(S\setminus A) \times P(\beta \setminus S)$. Suppose $G(A)$ is
$P(A)$-generic over $M$ and $G$ is
$P(S\setminus A)$-generic over $M[G(A)]$,
and $H$ is $P(\beta \setminus S)$-generic
over $M[G(A),G]$.
Then if $\tau \in M[G(A),G,H]$ has
strict level $S$, $\tau \nin M[G(A)]$.
So an object in a forcing extension has strict
level $S$ just in case there is a
term for $\tau$ in the forcing language
with strict level $S$,
and if $\tau$ is a term of strict level
$S$, then in any generic extension
the value of $\tau$ has strict level
$S$.  Consequently in our
construction we pass freely between
objects of strict level $S$ in a
generic extension and term of strict
level $S$ in the forcing language.

Many of the constraints required for
commutativity with morass maps are
expressed in terms of the strict level
of objects in a forcing extension (or
correspondingly, terms in the forcing
language).  For instance, in [\ref{Dumas}]
we define a term
function to be level if the
strict level of any term in the domain
equals the strict term of its image
under the function.  Such maps will
commute with morass maps in the manner
required by our construction.

\section{Constructing an $\R$-monomorphism
on a real-closed field} We wish to construct a set of terms in
the forcing language for adding
$\aleph_2$ generic reals that is forced
in all generic extensions to be an
$\R$-monomorphism on the
finite elements of an ultrapower of
$\R$. It is a result of B. Johnson
[\ref{Johnson}] that $\eta_1$-ordered
real-closed fields with cardinality
$\aleph_1$ are $\R$-isomorphic in
models of ZFC+CH. This result strengthens
the classical result that
$\aleph_1$-saturated real closed fields
of cardinality $\aleph_1$ are
isomorphic.  It is conceivable that in
a naive back-and-forth construction of
an order-isomorphism between
$\eta_1$-ordered real-closed fields
that a choice is made for an image (or
pre-image) of the function under
construction that precludes
satisfaction of $\R$-linearity.
\begin{definition}$($Full real-closed
field$)$\label{thm3.1} A real-closed
field $D$ is full iff for every finite
element, $r+\delta$, where $r\in \R$
and $\delta$ is infinitesimal, $r\in
D$.
\end{definition}
We will need to extend two results due
to B. Johnson [\ref{Johnson}] to meet
the requirement of $\R$-linearity in
the context of constructing term
functions using a morass.
\begin{lemma}$(B.
Johnson)$\label{thm3.2} Assume that $D$
and $I$ are full real-closed subfields
of $\eta_1$-ordered real-closed fields
$D^*$ and $I^*$ (resp.), $\phi:D\to I$
is an $\R$-monomorphism, and $r\in R$. Then
there is an extension of $\phi$,
$\phi^*$, that is an $R$-monomorphism of
the real closure of the field generated
by $D$ and $r$, $F(D,r)$, onto the real
closure of the field generated by $I$
and $r$, $F(I,r)$. Furthermore $F(D,r)$
(and consequently, $F(I,r)$) is full.
\end{lemma}
\begin{lemma}$(B.
Johnson)$\label{thm3.3} Let $D$, $D^*$,
$I$, $I^*$ and $\phi$ be as in Lemma
\ref{thm3.2}, $x\in D^*$ and assume
that the real closure of the field
generated by $D$ and $d$, $F(D,d)$, is
full. Let $y\in I^*$ be such that
\[ (\all d\in D)(d<x \iff \phi(d)<y). \]
Then there is an $\R$-monomorphism
extending $\phi$, $\phi^*:F(D,r)\to
I^*$, such that $\phi^*(x)=y$.
\end{lemma}
\begin{definition}$($Archimedean
valuation$)$\label{thm3.4} If $x$ and
$y$ are non-zero elements of a
real-closed field, they have the same
Archimedean valuation, $x \sim y$,
provided that there are $m, n \in
\mathbb{N}$ such that
\[ \mid x\mid <n\mid y\mid \]
and
\[ \mid y\mid<m\mid x\mid . \]
If $\mid x\mid <\mid y\mid$ and $x \nsim y$, then $x$
has Archimedean valuation greater than
$y$, $x \succeq y$.
\end{definition}
Archimedean valuation is an equivalence
relation on the non-zero elements of a
real-closed field (RCF). The non-zero
real numbers have the same valuation.
Elements with the same valuation as a
real number are said to have real
valuation. In a nonstandard real-closed
field, elements with valuation greater
than a real valuation are
infinitesimal.  The finite elements of
a real closed field are the
infinitesimal elements and those with
real valuation.

\section{The Esterle algebra}
We define the Esterle algebra
[\ref{Esterle2}] and review some basic
properties.
\begin{definition}
$(\soo)$\label{thm4.1} $\soo$ is the
lexicographic linear-ordering with
domain $\{s:\o1 \to 2 \mid s$ has
countable support and the support of s
has a largest element $\}$.
\end{definition}
\begin{definition}
$(\goo)$\label{thm4.2} $\goo$ is the
ordered group with domain $\{ g:\soo
\to \R \mid g$ has countable
well-ordered support$ \}$,
lexicographic ordering, and group
operation pointwise addition.
\end{definition}
We define an ordered algebra of formal
power series, $\ea$.  The universe of
$\ea$ is the set of formal power
series, $\sum_{\lambda <\gamma}
\alpha_{\lambda} x^{a_{\lambda}}$,
where:
\begin{enumerate}
\item $\gamma<\o1$.
\item $(\all \lambda<\gamma)\,
\alpha_{\lambda} \in \R$.
\item $\{a_{\lambda}\mid
\lambda<\gamma\}$ is a countable
well-ordered subset of $\goo$ and
$\lambda_1<\lambda_2<\gamma \Rightarrow
a_{\lambda_1}<a_{\lambda_2}$.
\end{enumerate}

The ordered algebra, $\ea$, is
isomorphic to the set of functions,
with countable well-ordered support,
from $\goo$ to $\R$.  The lexicographic
order linearly-orders $\ea$.  Addition
is pointwise and multiplication is
defined as follows:\, \\
Suppose
$a=\sum_{\lambda<\gamma_1}
\alpha_{\lambda} x^{a_{\lambda}}$ and
$b= \sum_{\kappa<\lambda_2}
\beta_{\kappa}x^{b_{\kappa}}$ are
members of $\ea$.  Let
\[ C= \{ c\mid (\some \lambda<\gamma_1)
(\some \kappa<\gamma_2) \, \,
c=a_{\lambda}+b_{\kappa} \}. \] Then
\[ a\cdot b= \sum_{c\in
C}((\sum_{a_{\lambda}+b_{\kappa}=c}
\alpha_{\lambda}\cdot \beta_{\kappa})
\, x^c ).
\]
\begin{definition}$($Esterle algebra,
$\ea)$\label{thm4.3} The Esterle
algebra, $\ea$, is $\{ f:\goo \to \R
\mid f$ has countable well-ordered
support$\}$. $\ea$ is lexicographically
ordered, with pointwise addition, and
multiplication defined above.
\end{definition}
$\R$ may be embedded in $\ea$ by
$\alpha \longmapsto \alpha x^e$, where
$e$ is the group identity in $\goo$.
Exponents in $\goo$ larger than $e$
(called positive exponents) correspond
to infinitesimal Archimedean
valuations, and those smaller than $e$
(called negative exponents) correspond
to infinite valuations.  The finite
elements of $\ea$ are those with
leading exponent $\geq e$.  The Archimedean valuations
of the Esterle algebra are represented by the group
of exponents  of $\ea$.
\begin{theorem}(J. Esterle [\ref{Esterle2}]) \label{thm4.4}
$\ea$ is an $\eta_1$-ordered
real-closed field.
\end{theorem}
A norm, $\| \, \|$, on an algebra $A$
is submultiplicative if for any $a,b
\in A$,
\[ \| a\cdot b\| \leq \| a\| \cdot \|
b\|. \]
\begin{theorem}(G. Dales [\ref{Dales1}], J. Esterle [\ref{Esterle3}]) \label{thm4.5}
The set of finite elements of $\ea$
bears a submultiplicative norm.
\end{theorem}
It is a standard result of model theory
that if $U$ is a non-principal
ultrafilter on $\omega$, the ultrapower
$\romu$ is an $\aleph_1$-saturated
real-closed field. Any two
$\aleph_1$-saturated, or
$\eta_1$-ordered, real-closed fields
with cardinality of the continuum are
isomorphic in models of ZFC+CH.  Hence
CH implies that there is a discontinuous homomorphism
of $C(X)$.
\begin{theorem}(B.Johnson [\ref{Johnson}]) \label{thm4.6}
(CH) If $U$ is a non-principal ultrafilter,
there is an $\R$-monomorphism
from the finite elements of $\romu$ into $\ea$.
\end{theorem}

We turn our attention to terms in a
forcing language $M^P$ that are forced
to be members of the Esterle algebra.
In [\ref{Dumas}] and [\ref{Dumas2}], we found sufficient
conditions for morass constructions.
The aggregate of these conditions were
characterized as morass-definability
and gap-2 morass-definability.
The central theorem of the papers were
that morass-definable
$\eta_1$-orderings are order-isomorphic
in the Cohen extension adding $\aleph_2$ generic reals
of a model of $ZFC + CH$
containing a simplified $(\o1,1)$-morass; and
gap-2 morass definable $\eta_1$-orderings
are order-isomorphic in the Cohen
extension adding $\aleph_3$ generic reals of
a model of $ZFC + CH$  containing a
simplified $(\omega_1,2)$-morass.
\begin{definition}(Level-dense)
Let $T_X\in M^P$ be forced to be a linear-ordering and
$X\subseteq M^P$ be a set of terms of strict level for the
domain of $T_X$.  $X$ is level-dense provided that for
$x, y\in X$, where the support of $x$ and the support of
$y$ are disjoint, and $G$, $P$-generic over $M$,
if $M[G]\models x<y$ then there is $z\in X\cap M$
such that $M[G] \models x<z<y$.
\end{definition}
\begin{definition}(Upward level-dense)
Let $T_X\in M^P$ be forced to be a linear-ordering,
and $X\subseteq M^P$ be a set of terms
of strict level for
the domain of $T_X$.  $X$ is upward level-dense provide that
for every $x,y,z\in X$, in which
$z$ has strict level $A\subseteq \kappa$, $p\in P$
with $p\forces x<z<y$ and $B\supseteq A$, there
is a term of strict level $B$, $w\in X$,
such that  $p\forces x<w<y$.
\end{definition}
\begin{definition}$($Morass-commutative$\,)$\label{def4.8}
Suppose $\morass$ is a simplified $(\omega_1,1)$-morass, $\lambda\leq\omega_1$
and $X\subseteq M^{P(\omega_1)}$.
We say that $X$ is morass-commutative beneath $\lambda$ if for any
$\zeta<\xi\leq\lambda$,
$f\in \mathcal{F}_{\zeta \xi}$ and
$x\in X\cap M^{P_{\zeta}}$, $f(x)\in X$.
We say that $X$ is morass-commutative if $X$ is morass-commutative beneath $\omega_1$.
\end{definition}
\begin{definition}$($Morass-definable$\,
)$\label{thm4.8} Let $\la X,<\ra \in
M[G]$ be a linear ordering. $X$ is
morass-definable if there is a set of terms
$T\subseteq M^P$
\begin{enumerate}
\item $T$ is a morass-commutative set of terms of strict
level, and $val_G(T)=X$
\item $T$ is level dense and upward
level dense
\item Every term of $T$ has countable
support (that is, any term in $T$ is
a term in $M^{P(A)}$, for some $A$
adding generic reals indexed by a
countable subset of $\omega_2$)
\end{enumerate}
\end{definition}
Let $E\subseteq M^P$ be the set of
terms of strict level for elements in the Esterle algebra
in the forcing language of the poset $P$.
It is routine to check that
$\ea$ is morass-commutative and that
every element of the Esterle algebra is
the interpretation of a term with
countable support.
\begin{lemma}
The Esterle algebra is level-dense.
\end{lemma}
Proof:
Let $a, b\in M^P$ be discerning terms for elements of the Esterle
algebra bearing disjoint
supports.  Let $G$
be $P$-generic over $M$
and $M[G]\models a<b$.
We wish to show that $a$ and $b$ are
separated in $M[G]$ by an element of $\ea\cap M$.
We work in $M[G]$.  Let $\gamma_1$ and $\gamma_2$ be
countable ordinals and
\[ a=\sum_{\lambda<\gamma_1}
\alpha_{\lambda} x^{a_\lambda} \]
and
\[ b=\sum_{\lambda<\gamma_2}
\beta_{\lambda} x^{b_\lambda}. \]
If $a$ and $b$ are equal on a partial sum, then that partial
sum is in $M$, so subtracting the largest common
partial sum of $a$ and $b$,
we may assume that $a$ and $b$ differ on
the first term of the sums, and
\[ \alpha_0 x^{a_0} \neq \beta_0
x^{b_0}. \]
If $a_0=b_0$, then $a_0 \in M$ and
there is $q\in \Q$ such that
\[ \alpha_0<q<\beta_0. \]
Then $q x^{a_0}\in M$ and
\[ a<q x^{a_0}<b. \]
Hence we assume that $a$ and $b$
have distinct Archimedean valuations.

We consider the case $a_0<b_0$.  Then $\alpha_0<0$. It is
sufficient to prove that there is an
element of $\ea \cap M$ that has valuation
between $a_0$ and $b_0$.  Let
\[ a_0=f:\soo \to \R \]
and
\[ b_0=g:\soo \to \R. \]
If $a_0$ and $b_0$ are equal on an
initial segment of their supports, then
this initial segment is in the ground
model.  We can therefore assume that
$a_0$ and $b_0$ either have distinct
least members, or have the same least
member of their supports, $s\in M$, and
\[ f(s)<g(s).  \]
In the latter case there is $q\in \Q$
such that
\[ f(s)<q<g(s). \]
Then $\{ (s,q) \} \in M$ and
\[ a_0< \{ (s,q) \}<b_0. \]
Hence we have left to consider the case
in which $s_{a_0}$ is the least member
of the support of $a_0$, $s_{b_0}$ is
the least member of the support of
$b_0$ and $s_{a_0}\neq s_{b_0}$.  If
either $s_{a_0}$ or $s_{b_0}$ are in
$M$ we can find a member of $\goo \cap
M$ that is a valuation between $a_0$
and $b_0$.  So we may assume that
neither $s_{a_0}$ nor $s_{b_0}$ are in
$M$. As we shall see, it is sufficient
to show that between any distinct
elements of $\soo$ from mutually
generic extensions, there is a member
of $\soo$ in the ground model.  We
assume without loss of generality that
$s_{b_0}<s_{a_0}$ (the case
$s_{a_0}<s_{b_0}$ is altogether
similar). Treating $s_{b_0}$ and
$s_{a_0}$ as countable subsets of
$\o1$, let $\mu$ be the least element
of $s_{b_0}$ that is not a member of
$s_{a_0}$. Let $\Delta=s_{a_0} \cap
\mu$.  Then $\Delta=s_{b_0}\cap \mu \in
M$. We note that $b_0(\mu)>0$,
otherwise $b_0<a_0$, contrary to
assumption. Let
\[ s_{c_0}=\Delta \cup \{ \mu \}. \]
Then $s_{c_0}\in M$ and
\[  s_{b_0}<s_{c_0}<s_{a_0} .\]
Let $c_0\in M$ be defined so that, for $\rho<\mu$,
\[ b_0(\rho)=c_0(\rho) \]
and
\[ 0<c_0(\mu)<b_0(\mu). \]
Then $c_0 \in M$ and
\[ a_0<c_0<b_0. \]
Since $a<b$, the coefficient of $a_0>0$.  Let
$q\in \Q$ and $0<q<a_0$.  Then $q x^{c_0} \in M$ and
\[ a<q x^{c_0}<b. \]
Therefore $\ea$ is level dense. \halmos
\begin{lemma}
The Esterle algebra is upward level-dense.
\end{lemma}
Let $A\subseteq B$ be countable subsets
of $\omega_1$.  Suppose $x, y, z\in
M^P$, with $x<z<y$, and $z$ has
level $A$.  It is sufficient to
show that there is an element $w\in
\ea$ of strict level $B$ such that
\[ x<w<y. \]
We may assume without loss of
generality that $z\in M$. If $x$ and
$y$ have an identical partial sum, then
$z$ must share that partial sum, and it
is therefore in $M$. We may subtract
the partial sum from all three formal
power series and pass to $x$ and $y$
that disagree on the first term of the
formal power series.  If $x$ and $y$
have the same Archimedean valuation,
$a$, then the first term of $z$ has
valuation $a$.  Let $\alpha\in \R$ be
of strict level $B$ and lie between
the initial coefficients of $x$ and
$y$. Let $w= \alpha x^a$.  Then $w$ has
strict level $B$ and
\[ x<w<y. \]
So we assume that $x$ and $y$ have distinct valuation.  If
$x$ and $z$ have the same valuation, $a$, and the
initial coefficient of $x$ is negative, let $\alpha \in \R$
be negative and greater than the the initial coefficient of $x$ and
have strict level $B$.
If the initial coefficient of $x$ is positive, let
$\alpha$ be positive and greater then the initial coefficient of $x$
and have strict level $B$.
In either case, let $w=\alpha x^a$.  Then
$w$ has strict level $B$ and
\[ x<w<y. \]
The cases for $x$ and $y$ having the same
valuation are similar.

If $x$ and $y$ have distinct valuation,
let $\alpha \in \R$ be positive
and have strict level $B$. Let
$w=\alpha \cdot z$.  Then $w$ has strict
level $B$ and
\[ x<w<y. \]
Therefore $\ea$ is upward level-dense. \halmos
\begin{theorem}
The Esterle algebra is morass-definable.
\end{theorem}
\section{Extendable functions}
In order to extend an
order-monomorphism by commutativity
with a splitting map we must
satisfy both algebraic conditions
and order constraints.  In the
next section we show that
the morass-commutative extension
of an $\R$-monomorphism,
satisfying certain technical
constraints (extendability), may be extended
to an $\R$-monomorphism
satisfying those
same constraints. The technical
constraints are those required for an
inductive construction along the
vertices of a simplified morasses.

We state the central technical
condition that permits the inductive
construction of the following sections.
We restrict our attention to functions from
terms for a standard ultrapower of $\R$
to terms for elements of the Esterle algebra.
\begin{definition} $($Extendable
Function$)$ \label{thm5.2}
Let $M$ be a c.t.m.of ZFC, $P$ be the poset adding generic
reals indexed by a countable ordinal
and $G$ be $P$-generic over $M$.
The term function,
$\phi:X\to Y$, is extendable provided that the following
are satisfied:
\begin{enumerate}
\item $X$ is a set of terms of strict level for a
subring of a standard ultrapower of $\R$
(over $\omega$) in $M[G]$ that has countable transcendence
degree over $\R$
\item $Y$ is a set of terms of strict level for elements of the
Esterle algebra that is closed under partial sums and contains coefficients
\item $\phi:X\to Y$ is a level term function that is forced to be an $\R$-isomorphism
(with respect to $\R \cap X$).
\end{enumerate}
\end{definition}
We interpret a splitting map on a
set of countable ordinals as a function
on terms in the forcing language as
in [\ref{Dumas}].
If $\nu$ is a vertex of a morass, $\theta_{\nu}$ is
the ordinal associated with $\nu$, $X\cap M^{P(\theta_{\nu})}$
is morass-commutative beneath $\nu$ and
$\sigma\in \mathcal{F}_{\nu \nu+1}$ is the splitting function on $\theta_{\nu}$
then $X\cup \sigma[X]$ is morass-commutative beneath $\nu+1$.
We will show that
an extendable function may be
commutatively extended by a splitting
function to an extendable function.
\section{Commutative extensions of term functions}
In the inductive construction of the
following sections we will need
technical lemmas of two types: those
insuring that commutativity with morass
maps may be used to extend extendable
functions to extendable functions,
and those allowing the
extension of the domain of an extendable function
by a specified element to an extendable function.
Throughout this section we assume:
\begin{enumerate}
\item $\phi:X\to Y$ is an extendable
function.
\item $\sigma$ is a splitting map with
splitting point $\delta=0$.  That is,
$\theta<\omega_1$, $\sigma:\theta \to
\omega_1$ and $\sigma[\theta]\cap
\theta =\emptyset$.
\item $P$ is the poset adding generic
reals indexed by $\theta$.
\item $G$ is $P$-generic over $M$.
\item $H$ is $\sigma[P]$-generic
over $M[G]$.
\end{enumerate}
We wish to show that the
morass-commutative extension of $\phi$
to the ring generated by $X \cup
\sigma[X]$ is extendable.
\subsection{Splitting maps and algebraic independence}
The central result of this subsection
states, roughly, that a subset of a
field in a generic extension that is
algebraically independent (AI) over the
restriction of the field to the ground
model, will be AI over the restriction
of an extension field to a mutually
generic extension. It will follow that
the union of an AI subset of a field in
a generic extension with its morass
``split" in a mutually generic
extension will be AI over the
restriction of the field to the ground
model.
\begin{lemma}\label{thm6.1}
Let $P$ be the poset adding
generic reals indexed by $\omega_1$ and
$G$ be $P$-generic
over $M$. Let $F$ be a morass-definable
real-closed field in $M[G]$.
Let $\theta<\omega_1$, $P(\theta)$ be the poset
adding generic reals indexed by
$\theta$, $G(\theta)$ be the $P(\theta)$-generic factor
of $G$, and $\sigma$ be a
splitting map on $\theta$ (with
splitting point $\delta=0$). If
$\chi=\{ x_1,\ldots,x_n\} \subseteq
(F\cap M[G])$ is linearly independent
(LI) over $M\cap F$ and $H$ is
$\sigma[P]$-generic over $M[G]$, then
$\chi$ is LI over $M[H]\cap F$.
\end{lemma}
Proof: Let
\[ \mathbf{x}=(x_1, \ldots, x_n)\in
F^n\cap M[G]. \] Suppose that there is
a non-zero vector of $\mathbf{y}\in
F^n\cap M[H]$ which is orthogonal to
$\mathbf{x}$.  Without loss of
generality we may assume that the
components of $\mathbf{y}$ are LI over
$M\cap F$. We consider terms for
$\mathbf{x}$ and $\mathbf{y}$, $x\in
M^P$ and $y\in M^{\sigma[P]}$
respectively. Let $p\in P$, $q\in
\sigma[P]$, be such that
\[ p\cdot q \forces x \cdot
y=0. \] We work below $p\cdot q$ and
assume that $x \cdot y=0$ in all
$P(\theta \cup \sigma[\theta])$-generic
extensions.  Let $\la \sigma_i \mid
i\in \omega \ra$ be a sequence of
splitting functions such that
\begin{enumerate}
\item $\sigma_0=\sigma$.
\item For any $i\in \omega$, $\sigma_i$
has splitting point $\delta=0$.
\item For any $i\in \omega$, the range
of $\sigma_i$, $\theta_i$, is contained
in $\omega_1$.
\item If $i\neq j$, then the range of
$\sigma_i$ is disjoint from the range
of $\sigma_j$.
\end{enumerate}
Let $P_i$ be the poset adding generic
reals indexed by $\theta_i$ and $Q_1$
be the poset adding generic reals
indexed by $\bigcup_{i\leq n}
\theta_i$. Let $q_1 \in Q_1$ be a
condition which forces that the span of
terms $\{ \sigma_i(y) \mid i \leq n \}$
has maximal possible dimension, $m\leq
n$, and let $\tau$ be a term for the
row-reduced echelon form of the matrix
with rows $\sigma_1(y), \ldots,
\sigma_n(y)$. We consider $q_1$ as a
condition in $\Pi_{i=1}^n P_i$,
$p_1\cdots p_n$, in which $p_i \leq q$
for all $i\leq n$. Hence $p\cdot q_1$
forces that $x$ is in the null space of
$\tau$.

For $0<i\leq n$, let
\[ \sigma_i^*=\sigma_{n+i} \comp
\sigma_i^{-1}. \] Then
\[ \sigma_i^*[\theta_i]=\theta_{i+n}. \]
Furthermore, through the definitions of
[\ref{Dumas}], $\sigma_i^*$ is
naturally interpreted as a function
from $P_i$ to $P_{i+n}$ and as well as
a function from $M^{P_i}$ to
$M^{P_{i+n}}$. Let
\[ \sigma^*:\bigcup_{i\leq n}
\theta_i \to \bigcup_{i\leq n}
\theta_{i+n}. \] If $Q_2$ is the poset
adding generic reals indexed by
$\bigcup_{i=1}^n \theta_{i+n}$ then
$\sigma^*:Q_1 \to Q_2$. Furthermore,
since $F$ is morass-definable,
$\sigma^*(\tau)$ is forced in all
$Q_2$-generic extensions to be a
row-reduced matrix with dimension $m$.

Let $H_1$ be $Q_1$ generic over $M$,
with $q_1\in H_1$ and $H_2$ be
$Q_2$-generic over $M[H_1]$, with
$\sigma^*(q_1)\in H_2$. Let $A_1$ be
the value of $\tau$ in $M[H_1]$ and
$A_2$ be the value of $\sigma(\tau)$ in
$M[H_2]$.  Then $A_1$ and $A_2$ are
row-reduced matrices with dimension
$m$.  Furthermore, since $m$ is the
maximum possible rank for $\tau$ in any
generic extension,
\[ Span(A_1)=Span(A_2). \]
The row reduced echelon form is
canonical, so
\[ A_1=A_2. \]
However $M[H_1]$ and $M[H_2]$ are
mutually generic extensions, $A_1\in
M[H_1]$ and $A_2\in M[H_2]$.  Therefore
\[ A_1=A_2\in M. \]
In $M[G,H_1,H_2]$, $x$ is in the
null-space of $A_1\in M$.  This
contradicts the assumption that the
components of $x$ are LI over $M\cap
F$. \halmos
\begin{lemma}\label{thm6.2}
Let $F$, $G$, $H$, $\chi$ satisfy the
hypotheses of Lemma \ref{thm6.1}.
Assume $n\in \N$,
$\bar{x}=(x_1,\ldots,x_n)\in F^n\cap
M[G]$, $\bar{y}=(y_1,\ldots,y_n)\in
F^n\cap M[H]$ and
\[ \sum_{i=1}^n x_i\cdot y_i = 0. \]
Then $\bar{y}$ is in the span of
vectors in $F^n\cap M$, all of which
are orthogonal to $\bar{x}$.
\end{lemma}
Proof:  Let \[ \bar{x}=(x_1,\ldots,x_n)
\]
and
\[ \bar{y}=(y_1,\ldots,y_n). \]
Similar to the previous argument, we
work beneath a condition that forces
$\bar{x} \perp \bar{y}$.  Let
$\sigma_1,\ldots,\sigma_m$ be splitting
functions with pairwise disjoint ranges
and
$\sigma_1(\bar{x}),\ldots\sigma_m(\bar{x})$
be a maximal LI set of vectors in
mutually generic extensions.  It is a
consequence of this maximal linear
independence that in any mutually
generic extension of $M$, the
interpretation of $\bar{x}$ is in the
span of $\sigma_1(\bar{x}),\ldots,
\sigma_m(\bar{x})$.  For $i\leq n$, let
\[ P_i=\sigma_i[P] \]
and \[ P^*=\Pi_{i=1}^m P_i. \] Let
$A_1$ be the rank $m$ matrix with row
vectors
$\sigma_1(\bar{x}),\ldots,\sigma_m(\bar{x})$
as the first $m$ rows, and the the zero
vector for the remaining rows.  The
nullity of $A_1$ is $n-m$.  Let $B_1$
be the row-reduced echelon form of the
basis matrix of the null space of
$A_1$. We note that, in a generic
extension including generic reals
indexed by $\sigma(P)$, $\bar{y}$ is in
the null space of $A_1$ and hence in
the span of the rows of $B_1$.  Let
$\tau_{A} \in M^{P^*}$ be a term for
$A_1$, and $\tau_{B}$ be a term for
$B_1$. We work beneath a condition of
$P^*$ that forces that $\tau_{B}$ is
the row-reduced echelon form of the
basis matrix of the null space of
$\tau_{A}$.  Let $G^*$ be $P^*$-generic
and $H^*$ be $P^*$-generic over
$M[G^*][H]$.  As in the previous Lemma,
if $A_2$ [resp. $B_2$] is the
interpretation of $\tau_A$ [resp.
$\tau_B$] in $M[H^*]$, then $A_2$, is a
row-reduced echelon matrix of rank $n$
and, in $M[H^*][H]$, $\bar{y}$ is in
the null space of $A_2$.  Additionally
$B_2$ is the row-reduced echelon form
of the basis matrix for the null space
of $A_2$. Every row of $A_2$ is in the
span of $A_1$, and both have dimension
$m$, so the range of $A_1$ and $A_2$
are identical.  Therefore
\[ B_1=B_2. \]
However $B_1\in M[G^*]$ and $B_2\in
M[H^*]$, with $G^*$ and $H^*$ mutually
generic.  Hence
\[ B_1 \in M. \]
Therefore $\bar{y}$ is in the span of
vectors in $F^n\cap M$. \halmos
\begin{corollary}\label{thm6.3}
Let $F$, $G$ and $H$ satisfy the
hypotheses of Lemma \ref{thm6.1}. If
$\chi \in M[G] \cap F$ is algebraically
independent (AI) over $M\cap F$, then
$\chi$ is AI over $M[H]\cap F$.
\end{corollary}
Proof: Let $\chi^*$ be the
multiplicative semi-group generated by
the elements of $\chi$. Then $\chi^*$
is LI over $M\cap F$. By Lemma
\ref{thm6.1}, $\chi^*$ is LI over
$M[H]\cap F$. Therefore $\chi$ is AI
over $M[H]\cap F$. \halmos
\begin{corollary}\label{thm6.4}
If $\chi$ is $AI$ over $M\cap F$, then
$\chi \cup \sigma[\chi]$ is AI over
$M\cap F$.
\end{corollary}
Proof: Let $G$ be $P$-generic over $M$,
and $H$ be $\sigma[P]$-generic over
$M[G]$.  $F$ is morass-commutative so
$\sigma[\chi]$ is AI over $M$. Suppose
there is a nontrivial linear
combination (over $M\cap F$) of
distinct elements of the semigroup
generated by $\chi \cup \sigma[\chi]$
that equals 0.  By Lemma \ref{thm6.4},
$\chi$ is AI over $M[H]\cap F$, so
there must be a nontrivial linear
combination (over $M\cap F$) of
elements of the semigroup generated by
$\sigma[\chi]$ that equals 0. However
$F$ is morass commutative, and
$\sigma[\chi]$ is AI over $M\cap F$.
\halmos
\begin{lemma}\label{thm6.6}
Let $X$ be a subring of a
standard ultrapower of $\R$, $\romu$.  Let
$X^*$ be the ring
generated by $X\cup \sigma[X]$,
where $\sigma$ is a splitting function.
If $\phi:X\to \ea$ is an extendable
$\R$-monomorphism, then
there is an $\R$-monomorphism,
$\phi^*:X^* \to \ea$, extending $\phi$
and $\sigma(\phi)$.
\end{lemma}
Proof: An element of $X^*$ may be
expressed as $\sum_{i=1}^n x_i\cdot
\sigma(y_i)$, for some $n\in \N$, and
$x_1,\ldots,x_n,y_1,\ldots,y_n \in X$.
Let $S$ be the set of expressions of
this form.  We define a function
$\psi:S\to \ea$ where
\[ \psi(\sum_{i=1}^n x_i\cdot
\sigma(y_i)) = \sum_{i=1}^n
\phi(x_i)\cdot \sigma(\phi(y_i)). \]
Let $\iota:S \to X^*$ be the natural
quotient map from the expressions of
$S$ to $X^*$.  The kernel of $\iota$ is
the set of expressions of $S$ that sum
to $0$ in $X^*$.  $\psi$ defines a ring
homomorphism on $X^*$ if and only if
for any $s$ in the kernel of $\iota$,
\[ \psi(s)=0. \]
For $i\leq n$ let $z_i=\sigma(y_i)$ and
$s=\sum_{i=1}^n x_i\cdot z_i$ be in the
kernel of $\iota$.  Then in $X^*$,
\[ \sum_{i=1}^n x_i\cdot z_i=0.
\]
Let
\[ \bar{x}=(x_1,\ldots,x_n)\in X^n \cap M[G], \]
\[ \bar{y}=(y_1,\ldots,y_n) \in X^n\cap
M[G] \] and \[ \bar{z}=(z_1,\ldots,z_n)
\in \sigma([X])^n \cap M[H]. \] By Lemma
\ref{thm6.2}, $\bar{z}$ is in the span
of elements of $X^n\cap M$ that are
orthogonal to $\bar{x}$.  Let $\{
b_1,\ldots,b_m\}$ be an LI set of
vectors of $X^n\cap M$ orthogonal to
$\bar{x}$ that contains $\bar{z}$ in
its span. Let $\la\cdot,\cdot \ra$ be
the dot product and
$(\alpha_1,\ldots,\alpha_m)\in \sigma([X])^n\cap
M[H]$ be such that \[ \sum_{i=1}^m
\alpha_i\cdot b_i=\bar{z}.
\] Let $\bar{\phi}:X^n\to \ea^n$ be defined by
\[
\bar{\phi}(s_1,\ldots,s_n)=(\phi(s_1),\ldots,\phi(s_n)).
\]
Recall that for a splitting map
$\sigma$, $\sigma(\phi):\sigma[X]\to
\ea$ is defined so that $\sigma$ and
$\phi$ commute. Then
\[ \psi(\la \bar{x},\bar{z}\ra) =
\la
\bar{\phi}(\bar{x}),\sigma(\bar{\phi})(\sum_{i=1}^m
\alpha \cdot b_i) \ra = \sum_{i=1}^m
\sigma(\alpha_i)\cdot \la
\bar{\phi}(\bar{x}),
\sigma(\bar{\phi})(b_i) \ra.
\] However $b_i\in X^n\cap M$ for all
$i\leq n$, so
\[ \sigma(\bar{\phi})(b_i)=\bar{\phi}(b_i). \]
Hence
\[ \sum_{i=1}^m
\sigma(\alpha_i)\cdot \la
\bar{\phi}(\bar{x}), \sigma(\bar{\phi})
(b_i) \ra = \sum_{i=1}^m
\sigma(\alpha_i)\cdot \la
\bar{\phi}(\bar{x}),\bar{\phi}(b_i) \ra
= \sum_{i=1}^m \sigma(\alpha_i)\cdot
\phi(\la \bar{x},b_i\ra). \] However,
for all $i\leq m$, $b_i\perp \bar{x}$.
So for all $i\leq m$,
\[ \phi(\la \bar{x},b_i\ra ) = 0
\]
and $\psi(\la \bar{x},\bar{z} \ra)=0$.
Therefore $\psi$ defines a ring
homomorphism on $X^*$, $\phi^*$, that
extends $\phi \cup \sigma(\phi)$ to
$X^*$.

Assume that $\bar{x}, \bar{y}\in
M[G]\cap X^n$,
$\bar{z}=\sigma(\bar{y})$ and  \[ \la
\bar{x},\bar{z} \ra \neq 0 \] then
there are \[ x_1',\ldots,x_m' \in X \]
and
\[ z_1',\ldots z_m' \in \sigma[X] \]
such that $\{ x_1',\ldots,x_m'\}$ is LI
over $M\cap X$ and
\[ \sum_{i=1}^m x_i'\cdot
z_i'=\sum_{i=1}^n x_i\cdot z_i \neq 0.
\] Then $\{
\phi(x_1'),\ldots,\phi(x_m')\}$ is LI
over $\phi[X] \cap M$.  By Lemma
\ref{thm6.1}, $\{
\phi(x_1'),\ldots,\phi(x_m')\}$ is LI
over $\ea \cap M[H]$.  Therefore
\[ \sum_{i=1}^n \phi^*(x_i\cdot z_i)=\sum_{i=1}^m
\phi^*(x_i'\cdot z_i')\neq 0. \] Thus
$\phi^*$ is a monomorphism.

Let $r\in \R \cap X^*$.  Then $r$ is an
element of the ring generated by
the reals of $\bar{X}
\cup \sigma[\bar{X}]$.
Hence, for all $r\in \R\cap X^*$,
\[ \phi^*(r)=r \]
and $\phi^*$ is $\R$-linear.
We show that $\phi^*$ is
order-preserving.
If $x_1,\ldots,x_n,y_1,\ldots,y_n \in X$,
we show that
$\sum_{i=1}^n x_i\cdot \sigma(y_i)>0$
iff $\sum_{i=1}^n \phi^*(x_i)\cdot
\phi^*(\sigma(y_i))>0$.

If $n=1$, then the sign of $x_1\cdot
\sigma(y_1)$ is the sign of the product
of the leading coefficients of $x_1$
and $\sigma(y_1)$, which are preserved
by $\phi$ and $\sigma(\phi)$
respectively. So
\[ x_1\cdot \sigma(y_1)>0 \iff \phi(x_1)\cdot
\sigma(\phi(y_1))>0.
\]
Assume that $n \geq 2$, and
\[ \sum_{i=1}^n x_i\cdot \sigma(y_i)>0.
\]
Let
\[ y=\phi^*(\sum_{i=1}^n x_i\cdot
\sigma(y_i))=\sum_{\lambda<\gamma}
\alpha_{\lambda}x^{a_{\lambda}}. \] By
assumption $\phi[X]$ is closed under partial
sums, so there are elements\\
$u_1,\ldots, u_j,u_{j+1},\ldots
u_k,v_1,\ldots,v_k \in Y$ such that
\begin{enumerate}
\item $y=\sum_{i=1}^k u_i\cdot
\sigma(v_i)$.
\item For $i\leq j$, every term in the
power series expansion of $u_i\cdot \sigma(v_i)$
has power (valuation) less than $a_0$.
\item For $j<i\leq k$, every term of
the power series expansion of $u_i\cdot \sigma(v_i)$
has power (valuation) at least $a_0$.
\end{enumerate}
Every term of $y$ expressed as a power series has
valuation no less than $a_0$, therefore
\[ \sum_{i=1}^j u_i\cdot \sigma(v_i)=0
\]
and
\[ y=\sum_{i=j+1}^k u_i\cdot
\sigma(v_i). \]
If $j<i\leq k$, then $s_i=u_i\cdot \sigma(v_i)$ has valuation
no less than $a_0$ in $\goo$.  If $s_i$ has valuation greater than $a_0$,
then $s_i$ has Archimedean valuation greater than $y$.  Let $S$ be the set
of indexes of the $u_i\cdot \sigma(v_i)$ with valuation $a_0$ and
$T$ be the set of indexes of the $u_i\cdot \sigma(v_i)$ with valuation
greater than $a_0$.  $S$ is nonempty and every element of
$T$ has Archimedean valuation greater than every element of $S$.

Since $\phi$ is extendable,
there are $b_0\in M[G]\cap \goo$ and $c_0\in M[H]\cap \goo$
such that $x^{b_0} \in \phi[X]$, $x^{c_0}\in \sigma[\phi[X]]$,
$b_0, c_0\geq e$ (in $G_{\omega_1}$) and
\[ a_0=b_0+c_0. \]
Let $u,v\in X$ be such that
\[ \phi(u)=x^{b_0} \]
and
\[ \sigma(\phi(v))=x^{c_0}. \]

We let $\bar{X}$ be the field generated by $X$
and $\psi:\bar{X} \to \ea$ be
the unique order-preserving
field-monomorphism extending $\phi$. We
have previously observed that $\phi^*$
is a ring monomorphism and $\psi \cup
\sigma(\psi)$ is an order-preserving
injection. So
\[ \phi^*(\sum_{i=1}^n x_i\cdot
\sigma(y_i))=\phi^*(\sum_{i=j+1}^k u_i\cdot \sigma(v_i)) \]
and
\[ \psi(\sum_{i=j+1}^k u_i/u\cdot
\sigma(v_i/v))=\alpha_0 + \sum_{0<\lambda<\gamma} \alpha_{\lambda}x^{(a_\lambda-a_0)}. \]
The real number $\alpha_0$ is in the domain and range of $\phi^*$,
$\phi^*$ is
$\R$-linear and $\phi^*\subseteq \psi$, so
\[
\phi^*(\alpha_0)=\psi(\alpha_0)=\alpha_0.
\]
Let $z=\sum_{0<\lambda<\gamma}
\alpha_{\lambda}x^{a_{\lambda}-a_0}$.

The range of $\phi^*$ is closed
under partial sums, so
$\sum_{0<\lambda<\gamma}
\alpha_{\lambda}x^a_{\lambda}$ is in
the range of $\phi^*$ and $u\cdot v$ is
in the range of $\psi$.  Thus
$z\in \psi[\bar{X}]$.
Every term of $z$ is infinitesimal.  Therefore
$\psi^{-1}(z)$ is infinitesimal and
$\sum_{i=1}^n x_i\cdot \sigma(y_i)$,
$\sum_{i=1}^n (x_i/u)\cdot \sigma (y_i/v)$
and $\alpha_0$ are each greater than $0$.
Therefore $\phi^*$ is extendable.  \halmos
\begin{lemma} \label{thm6.7} Let
\begin{enumerate}
\item $\bar{\nu}<\nu \leq \o1$
\item $D$ be a subring of a
standard ultrapower of $\R$
over $\omega$
\item $\phi:D\to \ea \in
M[G_{\bar{\nu}}]$ be an extendable
$\R$-monomorphism on $D$
\item $D^*=\bigcup_{\sigma \in
\mathcal{F}_{\bar{\nu}, \nu}}
\sigma[D]$.
\end{enumerate}
Then there is a unique extendable
$\R$-monomorphism $\phi^*$
on the ring generated by $D^*$ which,
for any $\sigma \in
\mathcal{F}_{\bar{\nu} \nu}$, extends
$\sigma[\phi]$.
\end{lemma}
Proof:  If $\nu=\bar{\nu}+1$, then the
result follows from Lemma \ref{thm6.6}.

If there is no limit ordinal ordinal
$\lambda$, $\bar{\nu}<\lambda \leq
\nu$.  Then there is $n\in \omega$ such
that
\[ \nu=\bar{\nu}+n. \]
By Lemma \ref{thm6.6}, for any
extendable $\phi:D\to \ea$ and
splitting function $\sigma$, the ring
monomorphism on the ring generated by
$D\cup \sigma[D]$ extending $\phi \cup
\sigma[\phi]$ is extendable. By $n$
iterated applications of Lemma
\ref{thm6.6}, there is a unique
extendable $\R$-monomorphic extension of $\phi$,
$\phi^*\supset \bigcup_{\sigma \in
\mathcal{F}_{\bar{\nu} \nu}}
\sigma[\phi]$, to the ring generated by
$\bigcup_{\sigma \in
\mathcal{F}_{\bar{\nu} \nu}}
\sigma[D]$.

So assume there is a limit ordinal
$\lambda$, $\nu_{\bar{\alpha}}<\lambda
\leq \nu_{\alpha}$. Let $\lambda$ be
the least limit ordinal greater than
$\nu_{\bar{\alpha}}$. Let
\[ D_{\lambda}=\bigcup_{\sigma \in
\mathcal{F}_{\bar{\nu}, \lambda}}
\sigma[D]. \] Let $D^*_{\lambda}$ be
the ring generated by $D_{\lambda}$, We
show that there is an extendable
$\R$-monomorphism of
$D^*_{\lambda}$ which, for any $\sigma
\in \mathcal{F}_{\bar{\nu} \lambda}$,
extends $\sigma[\phi]$.

Let $F$ be a finitely generated subring
of $D^*_{\lambda}$. Let $\{ d_1,\ldots,d_n\}$ generate $F$,
$\sigma_1, \ldots,\sigma_n \in
\mathcal{F}_{\bar{\nu},\lambda}$ and
for all $i\leq n$,  $c_i\in D$ be such
that
\[ d_i=\sigma_i(c_i). \]
By condition P4 in the definition of
the simplified morass, there is $n\in
\omega$, and $g\in
\mathcal{F}_{\nu_{\bar{\alpha}}+n,
\lambda}$ and $f_1,\ldots,f_n \in
\mathcal{F}_{\nu_{\bar{\alpha}},
\nu_{\bar{\alpha}}+n}$ such that, for
$i\leq n$,
\[ \sigma_i=g\comp f_i. \]
For each $m<n$, let $h_m$ be the
splitting function of
$\mathcal{F}_{\bar{\nu}+m,
\bar{\nu}+m+1}$. By Lemma \ref{thm6.6},
$\phi \cup h_1[\phi]$ may be extended to an
extendable $\R$-monomorphism. Furthermore this
ring monomorphism may be extended by
the splitting functions $h_2$ through
$h_m$.  Let $\psi$ be the function on
$\bigcup_{\sigma \in
\mathcal{F}_{\bar{\nu} \lambda}}
\sigma[D]$ resulting after the $n$
splits. Then $\psi$ is extendable and
$f_i(c_i)$ is in the domain of $\psi$
for all $i\leq n$. Therefore $g\comp
\psi$ is an $\R$-linear order
monomorphism and is the restriction of
$\phi^*$ to a ring containing $F$. Thus
\[ \phi^*_{\lambda}=\bigcup_{\sigma \in
\mathcal{F}_{\bar{\nu} \lambda}}
\sigma[\phi] \]
 is a well-defined extendable
$\R$-monomorphism of
$D^*_{\lambda}$.

By an inductive argument on $\nu$,
invoking condition P2, and the results
above at limits and Lemma \ref{thm6.6}
at successor ordinals, it is
straightforward to show that
$\bigcup_{\sigma \in
\mathcal{F}_{\lambda, \nu}}
\sigma[\phi^*_{\lambda}]$ has a unique
extension to an extendable $\R$-monomorphism
of the ring generated
by $\bigcup_{\sigma \in
\mathcal{F}_{\lambda, \nu}}
\sigma[D^*_{\lambda}]$. \halmos

\subsection{Extensions by a specified
element}

Because ordered subrings of real-closed
fields have unique extensions to
real-closed subfields, we will be able to
restrict our attention to extending
subrings by algebraically
independent elements. If
$X$ is a subring of a real-closed field
and $Z$ is a subset of that real-closed
field, we let $X[Z]$ be the subring
of the real-closed field
generated by $X\cup Z$.
We continue with conventions
similar to those outlined
earlier in the section.  Specifically
we assume:
\begin{enumerate}
\item $M$ is a c.t.m. of ZFC
\item $P$ is the poset adding generic
reals indexed by $\theta$
\item $G$ is $P$-generic over $M$
\item $F$ is the ring of finite elements of
a standard ultrapower of $\R$
over $\omega$ in M[G]
\item $X\subset F$
\item $Y\subset \ea$
\item $\phi:X\to Y$ is extendable.
\end{enumerate}
We wish to prove analogues of Johnson's
theorems that extendable
functions may be extended by a specified
element.
Throughout the arguments of this section, we will commonly
use $x$ to represent an element of $F$, and also
to represent the variable in the power series
representations of member of $\ea$.  Presumably the
context will make clear which use is intended.
\begin{lemma}\label{thm6.8}
Suppose $\phi:X\to Y$ is an extendable
function, $r\in \R$ and $r$ is transcendental
over $X$. Then there is an
extendable function that extends $\phi$, $\psi:X[r] \to
Y[r]$.
\end{lemma}
Proof: Let $\phi^*$ be
the $\R$-monomorphism
extending $\phi$ to the real closure of
$X$. Then $\phi^*\restrict_{\R\cap X}$ is $\R$-linear.
Suppose that $r$ is transcendental
over $X$.
The real closure of $X$ and the
real closure of $Y$ contain precisely the
same real numbers and are full.
By Lemma \ref{thm3.2},
there is an $\R$-monomorphic extension of
$\phi^*$, $\psi^*$, to the real closure
of the field generated by $X[r]$.  Let
$\psi=\psi^*\restrict_{X[r]}$.
Then $\psi$ is extendable. \halmos
\begin{lemma}\label{thm6.9}
Suppose $\phi:X\to Y$ is an extendable
function, $x\in F$ has strict level, and is transcendental over
$X$. Then there is $\bar{X}\supseteq
X[x]$ and extendable $\psi:\bar{X} \to \ea$
that extends $\phi$.
\end{lemma}
Proof: If
$x\in \R$, then apply
Lemma \ref{thm6.8}.
Assume $x\nin \R$.  We may assume that $x$ is infinitesimal.
Let $(l,u)$ be the gap formed by $x$ in
$X$, $L=\phi[l]$ and $U=\phi[u]$. The Esterle
algebra is an $\eta_1$-ordering, so
there is $y\in \ea$ that witnesses the
gap $(L,U)$. By application of
Johnson's Lemma \ref{thm3.3}, there is
$y\in \ea$ that witnesses the gap
$(L,U)$ and such
that the real closure of the field
extending $Y[y]$ is full. Although
the existence of such an element can be
used to advantage, the element $y$
may fail to have some of the properties
we require for a morass construction.

Let $\mu\subseteq \o1$ be the strict
level of $x$.
We seek a candidate in $\ea$ for the
image of $x$ with strict level $\mu$
that is transcendental over $Y$. The
Esterle algebra is level dense and
upward level dense, so by Lemma 4.5 of
[\ref{Dumas}] there is $y\in \ea$, with
the strict level $\mu$, such that for
all $z\in X$,
\[ x<z \iff y<\phi(z). \]
It may happen that $Y[y]$
is not closed under partial sums.\\ \\
Case 1: There is a largest partial sum of $y$
that is a member of $Y$.\\ \\
We include in this case that the first term of $y$
is a monomial not in $Y$.
Let $t$ be the largest partial sum
of $y$ that is also a member
of $Y$ and $s=\phi^{-1}(t)\in X$.
We shift our attention to the
gap formed by $x-s$ in $X$.
Then for all $z\in X$
\[ x-s<z \iff y-t<\phi(z). \]
We note that $x-s$ is transcendental over $X$.
Either $x-s$ has the same Archimedean valuation
as a member of $X$, or it has a distinct valuation.

Assume $x-s$ has the same Archimedean valuation
as a member of $X$. Let $\alpha x^a$ be the
leading term of $y-t$.
Then $x^a\in Y$ and $\alpha \nin Y$.
The coefficient $\alpha$ has strict level
contained in $\mu$.  Let $(L^*,U^*)$
be the cut formed in $\R \cap X$
by $\alpha$.  Then there is a
real number with strict level $\mu$
that witnesses the gap $(L^*,U^*)$.
So we assume that $\alpha$ has strict level $\mu$.
By Lemma \ref{thm6.8}
there is an extendable term function,
$\psi \supset \phi$, $\psi:X[\alpha]\to Y[\alpha]$,
with $\psi(\alpha \phi^{-1}(x^a))=\alpha x^a$.

We assume that $x-s$ has a
valuation distinct from the valuations
of members of $X$. We also note
that the strict level of $s$ and $t$
are equal and contained in $\mu$.
If the strict level of $s$ is a proper subset of
$\mu$, then $x-s$ and $y-t$ have level
$\mu$.  If the strict level of $x-s$ is
contained in $\mu$, then  there is an
element of $y^*\in \ea$ with a
valuation distinct from the valuations
of $Y$ having strict level equal to the
strict level of $x-s$.

We need consider only
the case in which $x$ has
valuation distinct from
the valuations of $X$.
Let $y\in \ea$ witness the gap
$(L,U)$, and
let $y_0=\alpha x^a$ be the leading
term of $y$.  If $\alpha>0$,
then $x^a$ witnesses
the gap $(L,U)$.

We assume without loss of generality
that $y_0=x^a$.
The strict level of $y_0$ equals
the strict level of $a\in \goo$.
It is straightforward to see that
there is $b\in \goo$ such that
\begin{enumerate}
\item $b$ is positive.
\item Any element of
the support of $b$ is greater
than any element of
the support of any exponent occurring in
any power series of $Y$.
\item $b$ has strict level $\mu$.
\end{enumerate}

Let $S$ be the union of supports of all
exponents among the elements in $Y$.
The exponent $b\in \goo$ is
greater than $0$ but less than any
positive exponent in $Y$.  Furthermore, $b$ is
greater than the constant $0$ function,
the additive identity of $\goo$ (and
the valuation of standard reals $\ea$).
However $b$ is greater than any positive
valuation occurring in $Y$. Consequently
$x^b$ is less than any positive infinitesimal
of $Y$. Additionally,
\begin{enumerate}
\item $x^{a+b}$ witnesses the gap
$(L,U)$.
\item $x^{a+b}$ is transcendental over
$X$.
\item $x^{a+b}$ has strict level $\mu$.
\item The ring generated by $Y\cup \{
x^{a+b}\}$ in $\ea$ is closed under partial sums and
coefficients.
\end{enumerate}
Let $\psi:X[x]\to Y[x^{a+b}]$ be the unique $\R$-monomorphism
extending $\phi$ such that
\[ \psi(x)=x^{a+b}. \]
Then $\psi$ is extendable.\\ \\
Case 2: There is no largest partial sum of $y$
that is a member of $Y$.\\ \\
Let $D$ be the countable well-ordered
sequence of exponents of $y$.  Let $D'$ be the smallest
initial segment of $D$ such that
$y\restrict_{D'}$ is not a member of
$Y$. Let $y'=y\restrict_{D'}$.  Then
$y'$ witnesses the gap
$(L,U)$. The strict level
of $y'$ is a subset of $\mu$, and
is possibly proper subset of $\mu$.  If the
strict level of $y'$ equals $\mu$, let
$\psi:X[x]\to Y[y']$ be the unique
$\R$-monomorphism
extending $\phi$ such that
\[ \psi(x)=y'. \]

Otherwise, let $a\in \goo$ be an
exponent of $\ea$ that has strict level
$\mu$ and is greater than all exponents
occurring in $Y$.  Let
\[ y^*=y'+x^{a}. \]
Then $y^*$ witnesses the gap $(L,U)$ and has
strict level $\mu$.  Let $\psi^*:X[x]\to Y[y^*]$
be the unique $\R$-monomorphism extending
extending $\phi$ such that
\[ \psi(x)=y^*. \]
Then $\psi$ satisfies the conditions for an
extendable function, except for closure under
partial sums. In particular $y'$ and $x^a$ are not
in the range of $\psi^*$.  Let $(L^*,U^*)$
be the gap formed by $x^a$ in $Y[y^*]$.  We observe
that $x^a$ is infinitesimal with respect to every
member of $Y[y^*]$.  Since $F$ is level dense and
upward level dense, there is a positive element of $F$,
$\epsilon$, with
strict level $\mu$ that is infinitesimal with respect
to all elements of $X[x]$.  Therefore there is
an $\R$-monomorphism
$\psi:X[x,\epsilon]\to Y[y^*, x^a]$ extending
$\psi^*$ and such that
\[ \psi(\epsilon)=x^a. \]
We note that $y\in Y[y^*,x^a]$ and
\[ Y[y^*, x^a]=Y[y, x^a]. \]
Then $Y[y, x^a]$ is closed
under coefficients and partial sums, so
$\psi$ is extendable. \halmos
These results permit a simplification of the construction.
Given an extendable $\R$-monomorphism, $\phi:D\to \ea$,
we may extend $\phi$ to an $\R$-monomorphism of $D[\R]$, and then
extend by algebraically independent infinitesimals.
\section{An $\R$-monomorphism from the finite elements of $\romu$
into the Esterle algebra.}
In the Cohen extension
adding $\aleph_2$ generic reals,
we construct a
level $\R$-monomorphism,
$\phi$, from the finite elements
of a morass-definable
ultrapower of $\R$ into $\ea$. This
construction differs significantly from
the construction of Woodin
[\ref{Woodin}]. The construction of
Woodin relies on the fact that in
the Cohen extension of $L$ by $\aleph_2$-generic
reals,
any cut of $\ea$ has a countable,
cofinal subcut.  As
Woodin observes, this argument is not
generalizable to models of ZFC with
higher powers of the continuum.  Our
construction yields a
monomorphism that is level, and therefore
respects the ``complexity" (with
respect to the index of Cohen reals)
of the elements of the ultrapower and
the Esterle algebra.  Consequently, for
any $S\subseteq \omega_2$, $\phi \cap
M[G(S)]$ is an $\R$-monomorphism
of ring of finite elements of $\romu \cap
M[G(S)]$ to $\ea \cap
M[G(S)]$.
\begin{theorem}\label{thm7.1}
Suppose $M$ is a c.t.m. of ZFC+CH
containing a simplified $(\o1,1)$-morass and $P$ is
the poset adding generic reals indexed
by ordinals less than $\omega_2$.  Let
$G$ be $P$-generic over $M$, $F\in
M[G]$ be the ring of finite elements of
an ultrapower of $\R$ over a standard
ultrafilter on $\omega$, and $\ea \in
M[G]$ be the Esterle Algebra computed
in $M[G]$. Then there is a level
$\R$-monomorphism,
$\phi:F\to \ea$.
\end{theorem}
Proof: For each $\nu$, $\alpha<\omega_1$,
let $G_{\nu}$ be the factor of $G$
adding generic reals indexed by
$\theta_{\nu}$ (the ordinal associated
with the vertex $\nu$ in the morass),
and, for $S\subseteq \omega_1$,
$G(S)$ be the factor of $G$
adding generic reals indexed by
$S$. Let
\[ F_{\nu}=F\cap M[G_{\nu}], \]
\[ F(S)=F\cap M[G(S)] \]
and $X$ be a maximal algebraically independent set of
positive morass generators of $F$.
Then $F$ is composed of
the finite elements of
the real closure of
the field generated by the morass-closure of $X$.

Let $\la x_{\alpha} \mid \alpha<\o1
\ra$ be a well-ordering of $X$. Let
$\la \nu_{\alpha} \mid \alpha<\omega_1
\ra$ be a weakly ascending transfinite
sequence of countable ordinals such
that $x_{\alpha} \in
M[G_{\nu_{\alpha}}]$. We will construct
an ascending sequence of functions
(ordered by inclusion), $\la
\phi_{\alpha}:D_{\alpha} \to \ea \mid \alpha<\omega_1
\ra$, such that for all $\alpha<\o1$
\begin{enumerate}
\item $x_\alpha \in D_{\alpha}$
\item $\phi_{\alpha}:D_{\alpha}\to
E_{\alpha}$ is extendable
\item If $x\in D_{\alpha}$,
$\nu \leq \nu_{\alpha}$, the morass
generator of $x$ is $x^*\in M[G_{\nu}]$
and $\sigma \in \mathcal{F}_{\nu
\nu_{\alpha}}$, then $\sigma(x^*)\in D_{\alpha}$.
\end{enumerate}
At each stage of the
construction, $\alpha<\o1$,
the domain of $\phi_{\alpha}$
extends the domain of
$\bigcup_{\beta<\alpha, f\in \script{F}_{\beta \alpha}}
f(\phi_{\beta})$, so that it
includes all morass descendants in
$M[G_{\nu_{\alpha}}]$,
and  $\phi_{\alpha}$
is extendable.

Case: $\alpha=0$.

Let $\la x_{0,n} \mid n\in \omega \ra
\subseteq M[G_{\nu_{0}}]$ be a
well-ordering of the morass-descendants
of $x_0$ in $M[G_{\nu_{0}}]$. We
construct a sequence of $\R$-monomorphisms, $\la
\phi_{0,n}:D_{0,n}\to \ea \mid n\in
\omega \ra$ such that for all $n\in
\omega$,
\begin{enumerate}
\item $D_{0,0}=\Q[x_{0,0}]$
\item $D_{0,n}[x_{0, n+1}]\subseteq D_{0, n+1}$
\item $\phi_{0,n}\in M[G_{\nu_0}]$
\item For all $m<n$, $\phi_{0,m}\subseteq
\phi_{0,n}$.
\item $\phi_{0,n}$ is extendable.
\end{enumerate}
Let $z=x_{0,0}$, and $\mu$ be the
strict level of $z$. We may assume that
$z$ is positive.  If $z\in \R$, let
$D_{0,0}=\Q[z]$ and
$\phi_{0,0}:D_{0,0}\to \ea$ be the
identity restricted to $D_{0,0}$.

If $z$ is infinitesimal, let $a\in
\goo$ be positive and have strict level
$\mu$, and
\[ y=x^a \in \ea. \]
Let $D_{0,0}=\Q[z]$ and
$\phi_{0,0}:D_{0,0} \to \ea$ be the
$\R$-linear ring monomorphism such that
\[ \phi_{0,0}(z)=y. \]

If $z=r+\epsilon$, where $\epsilon$ is
infinitesimal and positive, let $a\in
\goo$ have the same strict level as
$\epsilon$. Then let
\[ y=x^a \]
and $D_{0,0}=\Q[r,\epsilon]$. Let
$\phi_{0,0}:D_{0,0}\to \ea$ be the ring
monomorphism extending the identity on
$\Q[r]$ such that
\[ \phi_{0,0}(z)=y. \]
If $\epsilon<0$, let $y=-x^a$ and
define $D_{0,0}$ and $\phi_{0,0}$
analogously.

Let $N\in \omega$ and assume that $\la
\phi_{0,n} \mid n\leq N \ra$ satisfies
conditions 1 - 5 above (below $N+1$).
Then $\phi_{0,N}$ is extendable and
Lemmas \ref{thm6.8} and \ref{thm6.9}
apply. Let $z=x_{0,N+1}$.
Since $z$ is the image of a morass-generator
with strict level, $z$ is transcendental
over $D_{0,N}$.
If $z\in \R$, let
\[ D_{0,N+1}=D_{0,N}[z] \]
and
\[ \phi^*=\phi_{0,N} \cup \{
(z,z) \}. \] By Lemma \ref{thm6.8},
there is an $\R$-monomorphic extension of
$\phi^*$, $\phi_{0,N+1}:D_{0,N+1}\to
\ea$. The sub-ring of $\ea$ generated
by a set closed under partial sums and
coefficients is also closed under partial
sums and coefficients, so $\phi_{0,N+1}$
is extendable.

If $z$ is non-standard, let $R^*$ be
the set of reals contained in the
smallest full real closure of
$D_{0,N}\cup \{ z\}$.  By Lemma
\ref{thm6.8} there is an extendable
level $\R$-monomorphism
extending $\phi_{0,N}$,
$\phi^*:D_{0,N}[R^*]\to \ea$. Let
\[ D_{0,N+1}=D_{0,N}[R^*,z]. \] Then
$D_{0,N+1}\in M[G_{\nu_0}]$.  By Lemma
\ref{thm6.9} there is an extendable
extension of $\phi^*$,
\[ \phi_{0,N+1}:D_{0,N+1}\to \ea. \]
Let
\[ D_0=\bigcup_{n\in \omega} D_{0,n} \]
and
\[ \phi_0=\bigcup_{n\in \omega}
\phi_{0,n}. \]

Then the morass descendants of $x_0$
(in $M[G_{\nu_0}]$) are elements of
the domain of $\phi_0$, and $\phi_0$ is extendable.

Assume $\alpha$ a successor.

Let $\alpha=\bar{\alpha}+1$.  Assume
that $\la D_{\beta} \mid \beta < \alpha
\ra$ and $\la \phi_{\beta} \mid
\beta<\alpha \ra$ have been defined so
that for all $\gamma<\beta<\alpha$
\begin{enumerate}
\item $x_\gamma\in D_{\gamma}$
\item $D_{\gamma}\subseteq D_{\beta}$ is closed under morass-maps below $\beta$
\item $\phi_{\gamma} \subseteq
\phi_{\beta}$
\item $\phi_{\beta}:D_{\beta} \to \ea$ is an extendable $\R$-monomorphism.
\end{enumerate}
If $\nu_{\alpha}=\nu_{\bar{\alpha}}$,
then we may argue as in the previous
case.  Let $\la x_{\alpha,n} \mid
n<\omega \ra$ be an enumeration of the
morass descendants of $x_{\alpha}$ in
$M[G_{\nu_{\alpha}}]$.  We may extend
$D_{\bar{\alpha}}$ to $D_{\alpha}$
containing the morass descendants of
$x_{\alpha}$, and $\phi_{\bar{\alpha}}$
to $\phi_{\alpha}:D_{\alpha} \to \ea$
so that for all $\gamma<\beta \leq
\alpha$
\begin{enumerate}
\item $D_{\gamma} \subseteq D_{\beta}$.
\item $\phi_{\gamma}\subseteq
\phi_{\beta}$.
\item $\phi_{\alpha}$ is extendable.
\end{enumerate}
If $\nu_{\bar{\alpha}}<\nu_{\alpha}$,
then by Lemma \ref{thm6.7} there is an
extendable $\R$-monomorphism, $\phi^*$, extending
$\bigcup_{\sigma \in
\mathcal{F}_{\nu_{\bar{\alpha}}
\nu_{\alpha}}}
\sigma[\phi_{\bar{\alpha}}]$ to the
ring generated by $\bigcup_{\sigma
\in \mathcal{F}_{\nu_{\bar{\alpha}}
\nu_{\alpha}}}
\sigma[D_{\bar{\alpha}}]$, $D^*$. Let
$D'$ be ring generated by $D^*$ and
the morass descendants of $x_{\alpha}$
in $M[G_{\nu_{\alpha}}]$.  Let $\R^*$ be
the real numbers of the smallest full
extension of the real closure of $D'$.
Then the real closure of $D'[\R^*]$ is
full.  Let
\[ D_{\alpha}=D'[\R^*]. \]
By the preceding case, there is an
extendable $\R$-monomorphism
$\phi_{\alpha}:D_{\alpha} \to \ea$ with
$\phi_{\alpha}\supseteq \phi^*$. \\ \\
Finally, assume $\alpha<\o1$ is a limit
ordinal. \\ \\
If there is $\beta<\alpha$
such that $\nu_{\alpha}=\nu_{\beta}$
then we may proceed as in the case
$\alpha=0$ to define $D_{\alpha}$ and
$\phi_{\alpha}$.\\ \\
So we assume that
$\nu_{\beta}<\nu_{\alpha}$ for all
$\beta<\alpha$. Let
\[ \lambda=\bigcup_{\beta<\alpha}
\nu_{\beta}. \]
If $\lambda=\nu_{\alpha}$, then let
\[
D^*=\bigcup_{\beta<\alpha}(
\bigcup_{\sigma \in
\mathcal{F}_{\nu_\beta \nu_\alpha}}
\sigma[D_{\beta}]). \] Let $D$ be a
finitely generated subring of $D^*$
with generators $\{ d_1,\ldots,d_n\}$.
Then there is $\beta<\alpha$,
\[ C=\{ c_1,\ldots,c_n \} \]
and for, $i\leq n$, morass functions
$\sigma_i \in \mathcal{F}_{\nu_{\beta}
\nu_{\alpha}}$ such that
\[ \sigma_i(c_i)=d_i. \]
By condition P4 of Definition
\ref{thm2.1}, there is $\nu_{\beta}
\leq \gamma<\lambda$, $f_1, \ldots, f_n
\in \mathcal{F}_{\nu_{\beta} \gamma}$
and $g\in \mathcal{F}_{\gamma \lambda}$
such that for all $i\leq n$,
\[ \sigma_i(c_i)=g\comp f_i(c_i)=d_i.
\]
By Lemma \ref{thm6.7} there is an
extendable $\R$-monomorphism $\phi^*$ on the ring
generated by $\bigcup_{\sigma \in
\mathcal{F}_{\nu_{\beta} \lambda}}
\sigma[\phi_{\beta}]$.  It follows that
there is an extendable $\R$-monomorphism on the ring
generated by
$\bigcup_{\beta<\alpha}(\bigcup_{\sigma
\in \mathcal{F}_{\nu_{\beta}
\nu_{\alpha}}} \sigma[D_{\beta}])$
extending
$\bigcup_{\beta<\alpha}(\bigcup_{\sigma
\in \mathcal{F}_{\nu_{\beta}
\nu_{\alpha}}} \sigma[\phi_{\beta}])$.
We may then proceed as in earlier cases
to define an extendable $\R$-monomorphism
$\phi_{\alpha}:D_{\alpha} \to \ea$,
with $D_{\alpha}$ containing the morass
descendants of $x_{\alpha}$ in
$M[G_{\nu_{\alpha}}]$.

Finally, assume that
\[ \lim_{\beta<\alpha}
\nu_{\beta}=\lambda<\nu_{\alpha}. \] By
the previous argument, there is
extendable
$\phi^*:\bigcup_{\beta<\alpha}(\bigcup_{\sigma
\in \mathcal{F}_{\nu_{\beta} \lambda}}
\sigma[D_{\beta}]) \to \ea$ extending
the morass images in $M[G_{\lambda}]$
of the $\phi_{\beta}$. By Lemma
\ref{thm6.7} there is an extendable
$\R$-monomorphism
$\phi'_{\alpha}:D'_{\alpha} \to \ea$
such that
\[ D'_{\alpha} \supseteq
\bigcup_{\beta<\alpha}(\bigcup_{\sigma
\in \mathcal{F}_{\nu_{\beta}
\nu_{\alpha}}} \sigma[D_{\beta}]) \]
and
\[ \phi'_{\alpha} \supseteq \bigcup_{\beta<\alpha}(\bigcup_{\sigma
\in \mathcal{F}_{\nu_{\beta}
\nu_{\alpha}}} \sigma[\phi_{\beta}]).
\]
We proceed as in earlier cases to
extend $\phi'_{\alpha}$ to an
extendable $\R$-monomorphism,
$\phi_{\alpha}:D_{\alpha} \to \ea$,
where $D_{\alpha}$ extends $D'$ and
contains the morass descendants of
$x_{\alpha}$ in $M[G_{\nu_{\alpha}}]$.

Let $D^*$ be the ring in $M[G]$
generated by
$\bigcup_{\alpha<\o1}(\bigcup_{\sigma
\in \mathcal{F}_{\nu_{\alpha} \o1}}
\sigma[D_{\alpha}])$. Then by earlier
argument there is an $\R$-monomorphism, $\phi^*:D^*\to
\ea$, such that
\[ \phi^*\supseteq \bigcup_{\alpha<\o1}(\bigcup_{\sigma
\in \mathcal{F}_{\nu_{\alpha} \o1}}
\sigma[\phi_{\alpha}]). \]

Then $\phi^*$ is an $\R$-monomorphism
on a domain that contains
a transcendental basis for $F$, and
extends uniquely to an $\R$-monomorphism,
$\phi:F \to \ea$.
\halmos
\begin{theorem}
Suppose $M$ is a c.t.m.
of $ZFC+CH$ containing a simplified
$(\omega_1,1)$-morass and $M[G]$ is
the Cohen extensions adding $\aleph_2$
generic reals.  Then if $X$ is an
infinite compact Hausdorff space in $M[G]$,
there is a
discontinuous homomorphism of
$C(X)$ in $M[G]$.
\end{theorem}
Proof:  By Corollary 6.9 of
[\ref{Dumas}], any non-principal
ultrafilter on $\omega$ in $M$ may be
extended to a standard ultrafilter in
the Cohen extension adding
$\aleph_2$-generic reals.  If $U$ is a
standard ultrafilter, then there is a
level $\R$-monomorphism
from the finite elements of
$\romu$ into $\ea$.  Hence the finite
elements of $\romu$ bear a non-trivial
submultiplicative norm.  The theorem
follows from results of B. Johnson
[\ref{Johnson}].  \halmos

\section{Discontinuous homomorphisms of $C(X)$ in
a Cohen extension adding $\aleph_3$-generic reals.}

The existence of a discontinuous homomorphism of
$C(X)$ when the continuum has power $\aleph_2$ has been known for
a long while (Woodin [\ref{Woodin}]).  Unfortunately the methods
of the existing proof do not extend to higher powers of
the continuum.  It has been suggested that the
consistency of a discontinuous homomorphism of
$C(X)$ with higher powers of the continuum might
be proved with higher-gap morasses.   Theorem \ref{thm7.1}
is a modest strengthening of Woodin's result,
as the constructed $\R$-monomorphism is
an extendable function (except for the cardinality of the function).
We require a simplified $(\omega_1,2)$-morass for the
next construction (Velleman [\ref{Velleman2}]).
\begin{definition} $(Simplified \: \: (\kappa,2)-morass)$ \label{def3.1}
The structure $\gaptwo$ is a simplified $(\kappa,2)$-morass
provided it has the following properties:
\begin{enumerate}
\item $\gapone$ is a neat simplified $(\kappa^+,1)$-morass.
\item $\all \alpha<\beta\leq \kappa$,
$\script{F}_{\alpha \beta}$ is a family of embeddings (see page
172, $[\ref{Velleman2}]$) from $\langle \langle
\varphi_{\zeta}\mid \zeta<\theta_{\alpha}\rangle, \langle
\script{G}_{\zeta \xi} \mid \zeta<\xi\leq \theta_{\alpha}\rangle
\rangle$ to $\langle \langle \varphi_{\zeta}\mid
\zeta<\theta_{\beta}\rangle, \langle \script{G}_{\zeta \xi} \mid
\zeta<\xi\leq \theta_{\beta}\rangle \rangle$. \item $\all
\alpha<\beta<\kappa \: \: (\mid \script{F}_{\alpha \beta}
\mid<\kappa)$.
\item $\all \alpha<\beta<\gamma\leq \kappa \: \: (\script{F}_{\alpha \gamma}=\{ f \comp g\mid f\in \script{F}_{\beta \gamma}, g\in \script{F}_{\alpha \beta} \} )$. Here $f\comp g$ is defined by:\\
    \[ (f\comp g)_{\zeta}=f_{g(\zeta)} \comp g_{\zeta} \; \; \; \; for \; \zeta\leq \theta_{\alpha}, \]
    \[ (f\comp g)_{\zeta \xi}=f_{g(\zeta) g(\xi)} \comp g_{\zeta \xi} \; \; \; \; for \; \zeta<\xi\leq \theta_{\alpha}. \]
\item $\all \alpha<\kappa$, $\script{F}_{\alpha \alpha+1}$ is an
amalgamation (see page 173 $[\ref{Velleman2}]$). \item If
$\beta_1, \beta_2 <\alpha \leq \kappa$, $\alpha$ a limit ordinal,
$f_1\in \script{F}_{\beta_1 \alpha}$ and $f_2\in
\script{F}_{\beta_2, \alpha}$, then $\exists \beta (\beta_1,
\beta_2 < \beta <\alpha$ and $\exists f_1'\in \script{F}_{\beta_1
\beta} \: \exists f_2'\in \script{F}_{\beta_2 \beta} \: \exists
g\in \script{F}_{\gamma \alpha} (f_1=g\comp f_1'$ and $f_2=g\comp
f_2'))$. \item If $\alpha\leq \kappa$ and $\alpha$ is a limit
ordinal, then:
\begin{enumerate}
\item $\theta_{\alpha}=\bigcup \{ f[\theta_{\beta}] \mid
\beta<\alpha$, $f\in \script{F}_{\beta \alpha} \}$. \item $ \all
\zeta \leq \theta_{\alpha}$, $\varphi_{\zeta}=\bigcup \{
f_{\bar{\zeta }}[\varphi_{\bar{\zeta}}] \mid \exists \beta<\alpha
(f\in \script{F}_{\beta \alpha}$, $f(\bar{\zeta})=\zeta) \}$.
\item $\all \zeta<\xi \leq \theta_{\alpha}$, $\script{G}_{\zeta
\xi}=\bigcup \{ f_{\bar{\zeta} \bar{\xi}}[\script{G}_{\bar{\zeta}
\bar{\xi}}] \mid \exists \beta<\alpha \: (f\in \script{F}_{\beta
\alpha}$, $f(\bar{\zeta})=\zeta$, $f(\bar{\xi})=\xi ) \}$.
\end{enumerate}
\end{enumerate}
\end{definition}
\begin{theorem}\label{thm8.1}
Let $M$ be a c.t.m. of $ZFC+CH$
containing a simplified $(\omega_1,2)$-morass,
and $M[G]$ be a generic extension of $M$
adding $\aleph_3$ generic reals.
Let $X$ be an infinite compact Haussdorf space
in $M[G]$, and $C(X)$ be the algebra of continuous real-valued
functions of $X$ in $M[G]$.  Then there is a discontinuous homomorphism of $C(X)$
in $M[G]$.
\end{theorem}
Proof:
Let $M$ be a c.t.m. of $ZFC + CH$ containing
a simplified $(\omega_1,2)$-morass, $\gaptwo$.
Let $P$ be the poset adding generic reals
indexed by $\omega_3$, and $G$ be $P$-generic over $M$.
Then $\gapone$ is a simplified $(\omega_2,1)$-morass, and
below $\omega_1$, $\gapone$ satisfies the axioms of a
simplified $(\omega_1,1)$-morass.  Hence the construction
of Theorem \ref{thm7.1} below $\omega_1$
can be completed in  $M$.  In particular, if $U_0\in M$
is a non-principal ultrafilter in $M$, then by
Corollary 6.9 of [\ref{Dumas2}], there is $\overline{U}\subseteq M^{P(\omega_1)}$,
a standard term for an ultrafilter that is
morass-commutative below $\omega_1$,
that is forced to extend $U_0$.  Furthermore, the morass-continuation
of $\overline{U}$, $U$, is a standard ultrafilter that is
gap-2 morass-commutative and morass-commutative.
By Theorem 6.4 of [\ref{Dumas2}], $\romu$ is
a gap-2 morass-definable $\eta_1$ real-closed field.

We adapt the argument of Theorem \ref{thm7.1} to the gap-2 morass
construction of [\ref{Dumas2}].  Hence we construct a level term
function from the finite elements of a standard ultrapower
to the Esterle algebra that is closed under morass-embeddings and
is forced to be an $\R$-monomorphism.
For $\alpha<\omega_1$,
let $X_{\alpha}=(\romu) \cap M^{P_{\alpha}}$
and $Y_{\alpha}=\ea \cap M^{P_{\alpha}}$.
We consider $X_{\alpha}$ and $Y_{\alpha}$ as the restrictions of
$\romu$ and $\ea$, resp., to the forcing
language adding generic reals indexed by $\varphi_{\theta_{\alpha}}$.
In any $P$-generic extension of $M$, $M[G]$,
the interpretation of $X_{\alpha}$ in $M[G]$ is
the interpretation of $X_{\alpha}$ in $M[G_{\alpha}]$
where $G_{\alpha}$ is the factor of $G$ that is
$P_{\alpha}$-generic over $M$.

We construct a
morass-commutative level term injection from $X_{\omega_1}$
to $Y_{\omega_1}$ that is forced
to be an $\R$-linear order-monomorphism.  The closure
under embeddings, $f_{\theta_{\beta}}$ where $f\in \script{F}_{\beta \omega_1}$, of this
function will be the term function we seek.
Let $\{ x_{\beta} \mid \beta<\omega_1 \}\subseteq X_{\omega_1}$
be a transfinite sequence of terms
of strict level for a maximal algebraically independent set of
morass generators of the infinitesimal elements of $X_{\omega_1}$,
such that $x_{\alpha}\in X_{\alpha}$ for all $\alpha<\omega_1$.
We will inductively construct a transfinite sequence of
morass-commutative term functions $\langle F_{\beta}:D_{\beta}\to E_{\beta} \mid \beta<\omega_1 \rangle$
that satisfies the following for all $\alpha \leq \beta<\omega_1$,
\begin{enumerate}
\item $D_{\beta}\subseteq X_{\theta_{\beta}}$ is a morass-commutative subring
\item $E_{\beta}\subseteq Y_{\theta_{\beta}}$ is closed under partial sums and contains all coefficients
\item $D_{\alpha}\subseteq D_{\beta}$ and $E_{\alpha}\subseteq E_{\beta}$
\item $x_{\beta}\in D_{\beta}$
\item $F_{\beta}$ is a level term function that is forced to be an
order-preserving, $\R$-monomorphism
\item $f_{\theta_{\alpha}}[F_{\alpha}]\subseteq F_{\beta}$ for all
$f\in \script{F}_{\alpha \beta}$.
\end{enumerate}
We call a sequence of term functions satisfying these conditions (beneath $\beta$)
an extendable sequence.  We argue be induction on $\gamma<\omega_1$.

Base Case: $\gamma=0$.

Let $y_0$ be a positive infinitesimal monomial of $Y_0$ having the same strict level as $x_0$,
and $\R_0$ be the reals of the ground model.
Let $D_0$ be the ring generated by $\R_0 \cup \{ x_0\}$, $\R_0 [x_0]$,
and $E_0=\R_0 [y_0]$.  We observe that $E_0$ is closed under partial sums.
There is an $R$-linear order-monomorphism, $F_0:D_0 \to E_0$, with $F_0(x_0)=y_0$.

Successor Case: $\gamma=\beta+1$.  Let $\langle F_{\alpha}:D_{\alpha}\to E_{\alpha}
\mid \alpha \leq \beta \rangle$
be an extendable sequence satisfying conditions 1-6 above.
Let $D^*$ be the ring generated by
$\{ g_{\theta_{\beta}}[D_{\beta}] \mid g\in \script{G}_{\theta_{\beta} \theta_{\gamma}} \}$.
Then $D^*$ is generated by the union of the images of $D_{\beta}$ under the second components of
left-branching embeddings of $\script{F}_{\beta \gamma}$.
Let $h$ be the right-branching embedding of $\script{F}_{\beta \gamma}$  and
$D'$ be the ring generated by $D^*$ and $h_{\theta_{\beta}}[D_{\beta}]$.

By Lemma 5.2 of [\ref{Dumas2}],
$\bigcup \{ f_{\theta_{\beta}}[F_{\beta}] \mid f\in \script{F}_{\beta \gamma} \}$
is a level term injection that is forced to be an
$\R$-linear order-preserving injection.
The Lemma does not guarantee that
$\bigcup \{ f_{\theta_{\beta}}[F_{\beta}] \mid f\in \script{F}_{\beta \gamma} \}$
extends to a homomorphism of $D'$.  However, by Lemma \ref{thm6.7},
there is a unique extendable $\R$-monomorphism, $F^*:D^* \to \ea$.

Let $B$ be a transcendental basis of
$D_{\beta}$ over $\romu \cap M$
containing terms of strict level.
Let $D$ be the ring generated by $B$ and
$D_{\beta}\cap M$.  Let
$V$ be the semigroup generated by $B$.
Because $\script{F}_{\beta \gamma}$ is a set of
compatible embeddings, Lemma \ref{thm6.1} applies and,
treating $D$ as a vector space (over the field generated by
$D_{\beta} \cap M$), $V$ is a basis of elements
of strict level for $D$.
Therefore, if $F:D \to \ea$ is the naturally induced
$\R$-linear extension of
$\bigcup_{f\in \script{F}_{\beta \gamma}} f_{\theta_{\beta}}[F_{\beta}]$,
the image of $V$ under $F$ is a linearly independent
subset of $\ea$ over $E_{\beta} \cap M$.
Therefore $F$ is a linear transformation and
an $\R$-monomorphism of $D$.
We argue along the lines of Lemma \ref{thm6.6} to show that $F'$ is forced
to be order-preserving.
By Lemma \ref{thm6.6}, $F'\restrict_{D^*}$ is forced to be order-preserving.
If $z\in D'$, there is $n\in \N$,
$x_1,\ldots,x_n \in D^*$ and $y_1,\ldots,y_n \in f_{\theta_{\beta}}[D_{\beta}]$
such that $z=\sum_{i=1}^n x_i\cdot y_i$.
We show that
$\sum_{i=1}^n x_i\cdot y_i>0$
iff $\sum_{i=1}^n F'(x_i)\cdot
F'(y_i)>0$.

If $n=1$, then the sign of $x_1\cdot
y_1$ is the sign of the product
of the leading coefficients of $F'(x_1)$ and $F'(y_1)$.
However, $F'$ is $\R$-linear, so
\[ x_1\cdot y_1>0 \iff F'(x_1)\cdot
F'(y_1)>0.
\]
Assume that $n\geq 2$, and
\[ \sum_{i=1}^n x_i\cdot y_i>0.
\]
Let
\[ z=F'(\sum_{i=1}^n x_i\cdot
y_i)=\sum_{\lambda<\gamma}
\alpha_{\lambda}x^{a_{\lambda}}. \]
By assumption, the range of $F'$ is closed under partial
sums, so there are elements
$u_1,\ldots, u_j,u_{j+1},\ldots
u_k,v_1,\ldots,v_k \in Y$ such that
\begin{enumerate}
\item $z=\sum_{i=1}^k u_i\cdot
v_i$
\item For $i\leq j$, every term in the
power series expansion of $u_i\cdot v_i$
has power less than $a_0$
\item For $j<i\leq k$, every term of
the power series expansion of $u_i\cdot v_i$
has power at least $a_0$.
\end{enumerate}
Every term of $z$, expressed as a power series of $\ea$,
has valuation no less than $a_0$, therefore
\[ \sum_{i=1}^j u_i\cdot v_i=0
\]
and
\[ z=\sum_{i=j+1}^k u_i\cdot
v_i. \]
If $j<i\leq k$, then $s_i=u_i\cdot v_i$ has valuation
no less than $a_0$ in $\goo$.  If
the leading terms of $s_i$ has exponent greater than $a_0$,
then $s_i$ has Archimedean valuation less than $z$.  Let $S$ be the set
of all terms of all $s_i$, for all $1\leq j\leq i$, having exponent $a_0$ and
$T$ be the set of all terms of all $s_i$, $1\leq i\leq j$, with exponent
greater than $a_0$.  $S$ is nonempty and every element of
$T$ has Archimedean valuation less than every element of $S$.

There are $b_0, c_0\in \goo$
such that $x^{b_0}$ is in the range of $F^*$ and
$x^{c_0}$ is in the range of $f_{\theta_{\beta}}[F_{\beta}]$ and
\[ a_0=b_0+c_0. \]
Let $u\in D^*$ and $v\in f_{\theta_{\beta}}[D_{\beta}]$ be such that
\[ F'(u)=x^{b_0} \]
and
\[ F'(v)=x^{c_0}. \]

We observe that the ordered
ring $D'$ has a unique
extension to the field
closure of $D'$, $\bar{D}$.Treated as a field map,
the unique monomorphic extension of $F'$,
$\bar{F}:\bar{D} \to \ea$, is an $\R$-linear field monomorphism of $D'$.
We have previously observed that $F'$
is an $\R$-linear monomorphism, and $F^* \cup
f_{\theta_{\beta}}[F_{\beta}]$ is an order-preserving
injection. So
\[ F'(\sum_{i=1}^n x_i\cdot y_i)
=F'(\sum_{i=j+1}^k u_i\cdot v_i) \]
and
\[ \bar{F}(\sum_{i=j+1}^k u_i/u\cdot
v_i/v)=\alpha_0 + \sum_{0<\lambda<\gamma} \alpha_{\lambda}x^{(a_\lambda-a_0)}. \]
The real number, $\alpha_0$, is in the domain $\bar{F}$,
so $\bar{F}(\alpha_0)=F'(\alpha_0)=\alpha_0$.
The range of $F'$ is closed
under partial sums, so both
$\sum_{0<\lambda<\gamma}
\alpha_{\lambda}x^a_{\lambda}$
and $x^{a_0}$ are
in the range of $F'$.
Thus $\sum_{0<\lambda<\gamma}
\alpha_{\lambda}x^{a_{\lambda}-a_0}$ is
in the range of $\bar{F}$. Let
\[ \delta=\sum_{0<\lambda<\gamma}
\alpha_{\lambda}x^{a_{\lambda}-a_0}. \]
Every term of $\delta$ is infinitesimal (has Archimedean valuation
greater than elements of $\R$).
$\bar{F}$ is order preserving on $D'$,
so $\bar{F}^{-1}(\delta)$ is infinitesimal.
Hence
\[ \sum_{i=1}^n (x_i\cdot y_i)/(u\cdot v)=\alpha_0+\bar{F}^{-1}(\delta) \]
and the $\sum_{i=1}^n x_i\cdot y_i>0$ iff
$\alpha_0>0$ iff $F'(\sum_{i=1}^n x_i\cdot y_i)>0$.

Therefore $F':D' \to \ea$ is an $\R$-linear order-monomorphism, in which
the range of $F'$ contains all coefficients and all partial sums.
Let $\hat{D}$ be the ring generated by $D'$ and
$\R \cap M[G_{\theta_{\gamma}}]$.
By application of Lemma \ref{thm6.8}, we may extend
$F'$ to a level order-preserving $\R$-monomorphism,
$\hat{F}:\hat{D} \to \ea$.
Assume that $x_{\gamma}\nin D'[\R]$.
By Lemma \ref{thm6.9}, there is
$D_{\gamma}\supset \hat{D}$, with $x_{\gamma}\in D_{\gamma}$,
and an extension
of $\hat{F}$, $F_{\gamma}:D_{\gamma} \to \ea$,
such that $F_{\gamma}$ is an extendable function.

Limit Case: Suppose $\gamma$ is a limit ordinal.
Let $D'=\bigcup_{\beta<\gamma} (\bigcup_{f\in  \script{F}_{\beta \gamma}} { f_{\theta_{\beta}}[D_{\beta}]
})$. Then $D'$
is a subring of $D_{\omega_1}$ and has countable transcendence
degree over $\R$.  By condition 6 of Definition \ref{def3.1}
it is straightforward to verify that
$F':D' \to \ea$ defined by
$F'=\bigcup_{\beta<\gamma} (\bigcup_{f\in \script{F}_{\beta \gamma}} \{ f_{\theta_{\beta}}[F_{\beta}] \})$ is an extendable function.
Since $D'$ is full, we may apply Lemma \ref{thm6.8}
to extend $F'$ to
an $\R$-linear order-monomorphism of $\hat{D}$.
Let $D_{\gamma}=\hat{D}$.
Then by Lemma \ref{thm6.9} there is an extension of $F'$,
$F_{\gamma}:D_{\gamma}\to \ea$ that is an extendable function.
Let $F=\{ f_{\omega_1}[F_{\omega_1}]\mid f\in \script{F}_{\omega_1 \omega_2} \}$.
By Lemma 3.3 of [\ref{Dumas2}], the domain of $F$ is forced to be the finite elements of $\romu$.
\halmos
\begin{corollary}
Let $M$ be a c.t.m. of $ZFC+CH$
containing a simplified $(\omega_2,2)$-morass,
and $M[G]$ be a generic extension of $M$
adding $\aleph_4$ generic reals.
Let $X$ be an infinite compact Haussdorf space
in $M[G]$.  Then there is a discontinuous homomorphism of $C(X)$,
the algebra of continuous real-valued functions on $X$ in $M[G]$.
\end{corollary}
Proof.
By Woodin's argument [\ref{Woodin}],
in the generic extension adding $\aleph_2$
generic reals to $L$,
all cuts of $\ea$ have countable cofinal subcuts.
Let $M$ be a model of $ZFC+CH$ containing a
simplified $(\omega_2,2)$-morass
and $G$ be the generic extension adding
$\aleph_4$ generic reals.
If the cuts of $\ea$ in $M[G]$
admit countable cofinal subcuts, then
the argument for Theorem \ref{thm8.1} yields the existence of a
discontinuous homomorphism of $C(X)$
in a model with $2^{\aleph_0}=\aleph_4$.
\halmos

\section{Pressing the continuum}

The techniques of this paper, and those of [\ref{Dumas}] and [\ref{Dumas2}]
depend on the construction of functions
between sets of terms in
the forcing language of Cohen extensions,
utilizing commutativity with
order-preserving injections
indexing ordinals of the
Cohen poset.
Having presented the details of the construction
of level term functions using simplified gap-1 and
gap-2 morasses, it is relatively straightforward to
see how these constructions extend to
simplified higher (finite) gap morasses.

A simplified $(\kappa,n+1)$-morass, for
$\kappa$ a regular cardinal and integer
$n$, is a family of embeddings between
fake simplified $(\kappa, n)$-morass segments
that satisfies properties analogous to
those relating a simplified $(\omega_1,2)$-morass
to embeddings between fake simplified $(\omega_1,1)$-morass segments.
Central to the utility of these constructions is that
the second component of the embeddings
are order-preserving injections
between ordinals.  For a thorough
treatment of simplified finite gap morasses, including
an inductive definition, see Szalkai [\ref{Szalkai}].

Higher gap morasses will allow the extension of
results of this paper and [\ref{Dumas2}]
to Cohen extensions adding more than $\aleph_4$ generic reals.
The definition of gap-2 morass-definable $\eta_1$-orderings
and $\eta_1$-ordered real-closed fields (resp.) are easily generalized
to gap-n morass-definable $\eta_1$-orderings and
$\eta_1$-ordered real-closed fields (resp.).
We state the following without proof.
\begin{claim}
Let $M$ be a model of $ZFC + CH$ containing a
simplified $(\omega_1,n)$-morass, $P$ be the poset
adding $\aleph_n$ generic reals and $G$ be $P$-generic
over $M$.  Then in $M[G]$
\begin{enumerate}
\item Gap-n morass-definable $\eta_1$-orderings without endpoints
are order-isomorphic.
\item There exists a gap-n morass-definable, non-principal ultrafilter
over $\omega$.
\item If $U$ is a gap-n morass-definable non-principal ultrafilter over
$\omega$, then $\romu$ is a gap-n morass-definable $\eta_1$-ordered
real-closed field.
\item There exists an $\R$-isomorphism
between gap-n morass-definable $\eta_1$-ordered real-closed fields.
\item There exists an $\R$-monomorphism from
the finite elements of a standard ultrapower of $\R$ (over $\omega$)
into the Esterle algebra.
\item If $X$ is an infinite compact Hausdorff space, then there exists
a discontinuous homomorphism of $C(X)$, the algebra of continuous
real-valued functions on $X$.
\end{enumerate}
\end{claim}


\begin{thebibliography}{99}

\bibitem{Dales}\label{Dales1}
H.G. Dales, ``A discontinuous
Homomorphism of C(X)", Amer. J. Math.
Vol. 101, No. 3 (Jun., 1979), pp. 647-734

\bibitem{Dales2}\label{Dales2}
H.G. Dales, ``Automatic Continuity: A
Survey", Bull. Lond. Math. Soc. 10
(1978), 129 - 183.

\bibitem{Dumas}\label{Dumas}
B. A. Dumas, ``Order-isomorphic
$\eta_1$-orderings in Cohen
extensions", Annals of Pure and Applied
Logic, 158 (2009) 1 - 22.

\bibitem{Dumas2}\label{Dumas2}
B. A. Dumas, ``Morass-definable
$\eta_1$-orderings in Cohen extensions",
arXive 1701.02031 (2017).

\bibitem{Erdos}\label{Erdos}
P. Erd$\ddot{o}$s, L. Gillman and M.
Henriksen, ``An isomorphism theorem for
real closed fields", Ann. of Math. 61
(1955), 552 - 554

\bibitem{Esterle1}\label{Esterle1}
J. Esterle, ``Semi-normes on $C(K)$",
Proc. London Math. Soc. (3), 36 (1978),
27 - 45.

\bibitem{Esterle2}\label{Esterle2}
J. Esterle, ``Sur
l'existence d'un homomorphisme
discontinu do C(K)", Acta Math. Acad.
Sci Hungar., 30 (1977), 113 - 127.

\bibitem{Esterle3}\label{Esterle3}
J. Esterle, ``Injections de
semi-groupes divisible dans des
alg$\acute{e}$bras de convolution et
construction d'homomorphismes
discontinus des C(K)", Proc. London
Math. Soc. (3) 36 (1978), 59 - 85.

\bibitem{Hahn}\label{Hahn}
H. Hahn, ``$\ddot{U}$ber die
nichtarchimedishen
Gr$\ddot{o}$ssensysteme", S.B. Akad.
Wiss. Wien II(a) 116 (1907) 37 - 58.

\bibitem{Jech}\label{Jech}
T. Jech, \emph{Set Theory},
Springer-Verlag 2002.

\bibitem{Johnson}\label{Johnson}
B. E. Johnson, ``Norming $C(\omega)$
and related algebras", Trans. Amer.
Math. Soc., 220 (1976), 37 - 58.

\bibitem{Maclane}\label{Maclane}
S. Maclane, ``The universality of power
series fields", Bull. Amer. Math. Soc.
45 (1939) 888 - 890.

\bibitem{Szalkai}\label{Szalkai}
I. Szalkai, ``An Inductive Definition of
Higher Gap Simplified Morasses",
Publicationes Mathematicae Debrecen,
58/4 (2001), 605-634.

\bibitem(Velleman)\label{Velleman}
D. Velleman, ``Simplified Morasses",
Journ.  of Sym. Log., 49 (1984), 257 -
271.

\bibitem{Velleman2}\label{Velleman2}
D. Velleman, ``Simplified Gap-2 Morasses",
Annals of Pure and Applied Logic,
34 (1987), 171 - 208.

\bibitem{Woodin}\label{Woodin}
H. Woodin, ``A Discontinuous Homomorphism from C(X) without CH"
J. London Math. Soc. (1993) s2-48 (2): 299-315

\end{thebibliography}
\end{document}